\documentclass[preprint]{elsarticle}
\usepackage{mathrsfs}
\usepackage{amsmath}
\usepackage{amsfonts}
\usepackage{amssymb}
\usepackage[hidelinks]{hyperref}
\usepackage{subfigure}
\usepackage{color}
\usepackage{fancyhdr}
\usepackage{bbm}
\usepackage{bm}
\usepackage{graphicx}
\journal{Journal of Number Theory}
\newtheorem{definition}{Definition}[section]
\newtheorem{theorem}{Theorem}[section]
\newtheorem{remark}{Remark}[section]
\newtheorem{cor}{Corollary}[section]
\newtheorem{lemma}{Lemma}[section]
\numberwithin{equation}{section}

\begin{document}

\begin{frontmatter}
\title{A Weighted Divisor Problem\tnoteref{mytitlenote}}
\tnotetext[mytitlenote]{The work of Lirui Jia is supported by the National Natural Science Foundation of China (Grant No. 11571303). Wenguang Zhai is  supported by the National Key Basic Research Program of China (Grant No. 2013CB834201), the National Natural Science Foundation of China (Grant No. 11171344).}

\author{Lirui Jia}
\address{School of Mathematical Sciences, Zhejiang University,\\
Hangzhou 310027, People's Republic of China }

\ead{jialirui@126.com}

\author{Wenguang Zhai}
\address{Department of Mathematics, China University of Mining and Thechnology,\\ Beijing 100083, People's Republic of China }

\ead{zhaiwg@hotmail.com}

\date{}{}

\begin{abstract}
We study a weighted divisor function $\mathop{{\sum}'}\limits_{mn\leq x}\cos(2\pi m\theta_1)\sin(2\pi n\theta_2)$, where $\theta_i (0<\theta_i<1)$ is a rational number. By connecting it with the divisor problem with congruence conditions, we establish an upper bound, mean-value, mean-square and some power-moments.
\end{abstract}

\end{frontmatter}


\section{Introduction and main results}

\subsection{Introduction}

Let $d(n)=\sum\limits_{n=n_1n_2}1$ denote the divisor function, and $D(x)=\mathop{{\sum}'}\limits_{n\leq x}d(n)=\mathop{{\sum}'}\limits_{n_1n_2\leq x}1$ be the summatory function, where the prime $'$ on the summation sign indicates that if $x$ is an integer, then only $\frac{1}{2}d(x)$ or $\frac{1}{2}$ for $n_1n_2=x$ is
counted. In 1849, Dirichlet first proved that
$$D(x)=x\log x+(2\gamma-1)x+O(\sqrt{x}), \forall x\geq1,$$
where $\gamma$ is the Euler constant.

 Let
$$\Delta(x)=D(x)-x\log x-(2\gamma-1)x-\frac{1}{4}$$ be the error term in the asymptotic formula for $D(x)$. Dirichlet's divisor problem consists of determining the smallest $\alpha$, for which $\Delta(x)\ll x^{\alpha+\varepsilon}$ holds for any $\varepsilon>0$. Clearly, Dirichlet's result above implies that $\alpha\leq\frac{1}{2}$. Throughout the past more that 160 years, there have been many improvements on this estimate.
 The best estimate to-date has been given by
Huxley\cite{Huxley03,huxley2005exponential}, and reads
\begin{equation}\label{huxley}
\Delta(x)\ll x^\frac{131}{416}\log^\frac{26947}{8320}x.
\end{equation}
It is widely conjectured that α$\alpha=\frac{1}{4}$ is admissible and the best possible.

Since $\Delta(x)$ exhibits considerable fluctuations, one natural way to study the upper bounds is to consider the moments.

Voronoi's work\cite{voronoi1904fonction} in 1904 showed that
\begin{equation*}
  \int_1^X\Delta(x)dx=\frac{X}{4}+O(X^{\frac{3}{4}}).
\end{equation*}
Later, in 1922 Cram\'er\cite{cramer1922zwei} proved the mean square formula
\begin{equation*}
   \int_1^X\Delta(x)^2dx=cX^{\frac{3}{2}}+O(X^{\frac{5}{4}+\varepsilon}),\quad\forall~\varepsilon>0,
\end{equation*}
where $c$ is a positive constant. In 1983, Ivic \cite{ivic1983large} used the method of large values to prove that
\begin{equation}\label{ivic}
    \int_1^X|\Delta(x)|^Adx\ll X^{1+\frac{A}{4}+\varepsilon},\quad\forall~\varepsilon>0
\end{equation}
for each fixed $0\leq A\leq\frac{35}{4}$. The range of $A$ can be extended to $\frac{262}{27}$ by the estimate \eqref{huxley}.
In 1992, Tsang\cite{tsang1992higher}  obtained the asymptotic formula
\begin{equation}\label{01}
  \int_1^X\Delta(x)^kdx=c_kX^{1+\frac{k}{4}}+O(X^{1+\frac{k}{4}-\delta_k}),\quad \text{for}\ k=3,4,
\end{equation}
with positive constants $c_3$, $c_4$, and $\delta_3=\frac{1}{14}$, $\delta_4=\frac{1}{23}$. Ivi\'c and Sargos~\cite{ivic2007higher} improved the values $\delta_3$, $\delta_4$ to $\delta'_3=\frac{7}{20}$, $\delta'_4=\frac{1}{12}$, respectively.
Heath-Brown\cite{heath1992distribution} in 1992 proved that for any integer $k<A$, where $A$ satisfies \eqref{ivic}, the limit
$$c_k=\lim_{X\rightarrow\infty}X^{-1-\frac{k}{4}}\int_1^X\Delta(x)^kdx$$
exists.  Then, there followed a series
of investigations on explicit asymptotic formulas of the type \eqref{01} for larger values of $k$. In 2004, Zhai \cite{zhai04} established asymptotic formulas for $3\leq k\leq9$.

 At the beginning of 20th  century, Vorono\"{\i}\cite{voronoi1904fonction} proved the remarkable exact formula that
  \begin{equation*}
    \Delta(x)=-\frac{2}{\pi}\sqrt{x}\sum_{n=1}^\infty\frac{d(n)}{\sqrt{n}}\big(K_1(4\pi\sqrt{nx})+ \frac{\pi}{2}Y_1(4\pi\sqrt{nx})\big),
  \end{equation*}
  where $K_1$ is the modified Bessel function of the second kind of order 1, and $Y_1$ is the Bessel function of the second kind of order 1. The series on the right-hand side is boundedly convergent for $x$ lying in each fixed closed interval.

Recently, Berndt et al \cite{berndt2013weighted} considered $\mathop{{\sum}'}_{mn\leq x}\cos(2\pi m\theta_1)\sin(2\pi n\theta_2)$, which can be seen as a weighted divisor function and they got an analogue of Voronoi's formula as follows.

Let $J_1$ be the ordinary Bessel function. If $0<\theta_1, \theta_2<1$ and $x>0$,  then
\begin{align*}
   &\mathop{{\sum}'}\limits_{mn\leq x}\cos(2\pi m\theta_1)\sin(2\pi n\theta_2) \\
   =&  -\frac{\cot(\pi\theta_2)}{4}+\frac{\sqrt{x}}{4}\sum_{m=0}^\infty\sum_{n=0}^\infty \bigg\{\frac{J_1\big(4\pi\sqrt{(m+\theta_1)(n+\theta_2)x}\big)}{\sqrt{(m+\theta_1)(n+\theta_2)}}\\
   &+\frac{J_1\big(4\pi\sqrt{(m+1-\theta_1)(n+\theta_2)x}\big)}{\sqrt{(m+1-\theta_1)(n+\theta_2)}} -\frac{J_1\big(4\pi\sqrt{(m+\theta_1)(n+1-\theta_2)x}\big)}{\sqrt{(m+\theta_1)(n+1-\theta_2)}}\\ &-\frac{J_1\big(4\pi\sqrt{(m+1-\theta_1)(n+1-\theta_2)x}\big)}{\sqrt{(m+1-\theta_1)(n+1-\theta_2)}}\bigg\}.
\end{align*}

In this paper, we will study $\mathop{{\sum}'}_{mn\leq x}\cos(2\pi m\theta_1)\sin(2\pi n\theta_2)$ further with a different approach and prove some results analogous to those for $\Delta(x)$.

\textsc{Notations}. For a real number $t$, let $[t]$ be the largest integer no greater than $t$, $\{t\}=t-[t]$, $\psi(t)=\{t\}-\frac{1}{2}$, $\|t\|=\min(\{t\}$, $1-\{t\})$, $e(t)=e^{2\pi it}$. $\mathbb{C}$, $\mathbb{R}$, $\mathbb{Z}$, $\mathbb{N}$ denote the set of complex numbers, of real numbers, of integers, and of natural numbers, respectively; $C^r[a,b]$ be the class of functions having a continuous $r$th derivative in $[a,b]$; $f\asymp g$ means that both $f\ll g$ and $f\gg g$ hold. Throughout this paper, $\varepsilon$, $\varepsilon'$ denote sufficiently small positive constants, and $\mathcal{L}$ denotes $\log T$.

\subsection{Main results}\label{theorems}

We first state the divisor problem with congruence
conditions.  A divisor function with congruence conditions is defined by
\begin{equation*}
  d(n; r_1, q_1, r_2, q_2) =\mathop{{\sum}'}_{\begin{subarray}{c}  n=n_1n_2\\n_i\equiv r_i\!\!\!\!\!\pmod{q_i}\\i=1, 2\end{subarray}}{1},
  \end{equation*}
 of which, the summatory function is
 \begin{equation*}
  D(x; r_1, q_1, r_2, q_2) =\mathop{{\sum}'}_{\begin{subarray}{c}  n_1n_2\leq x\\n_i\equiv r_i\!\!\!\!\!\pmod{q_i}\\i=1, 2\end{subarray}}{1}.
\end{equation*}

From~Richert \cite{Richert}, we can find that for $x\geq q_1q_2$, $1\leq r_i\leq q_i$ ($i=1, 2$)
\begin{align}\label{4}
   &D(x; r_1, q_1, r_2, q_2)\\
   \nonumber=& \frac{x}{q_1q_2}\log\Big(\frac{x}{q_1q_2}\Big) - \bigg(\frac{\Gamma'}{\Gamma}\Big(\frac{r_1}{q_1}\Big) + \frac{\Gamma'}{\Gamma}\Big(\frac{r_2}{q_2}\Big) +1\bigg)\frac{x}{q_1q_2}+\Delta(x; r_1, q_1, r_2, q_2),
\end{align}
where
\begin{equation*}
    \Delta(x; r_1, q_1, r_2, q_2)=-\!\!\!\!\sum_{\begin{subarray}{c}n_1\leq \sqrt{\frac{q_1 x}{q_2}}\\n_1\equiv r_1\!\!\!\!\!\pmod{\!q_1}\end{subarray}}\!\psi\big(\frac{x}{q_2n_1}-\frac{r_2}{q_2}\big)-\!\!\!\!\sum_{\begin{subarray}{c}n_2\leq \sqrt{\frac{q_2x}{q_1}}\\n_2\equiv r_2\!\!\!\!\!\pmod{\!q_2}\end{subarray}}\!\psi\big(\frac{x}{q_1n_2}-\frac{r_1}{q_1}\big)+O(1).
\end{equation*}

Denote $S(x; \theta_1, \theta_2)=\mathop{{\sum}'}_{mn\leq x}\cos(2\pi m\theta_1)\sin(2\pi n\theta_2)$. Let $\theta_1$, $\theta_2$ be rational numbers satisfying $\theta_i=\frac{\displaystyle a_i}{\displaystyle q_i}$, with $(q_i, a_i)=1$ ($i=1, 2$). Then $S(x; \theta_1, \theta_2)=S(x; \frac{\displaystyle a_1}{\displaystyle q_1}, \frac{\displaystyle a_2}{\displaystyle q_2})$. For simplicity, in this paper we consider $S(q_1q_2x; \frac{\displaystyle a_1}{\displaystyle q_1}, \frac{\displaystyle a_2}{\displaystyle q_2})$.

Our first result is an analogue of \eqref{ivic} for $\Delta(x; r_1, q_1, r_2, q_2)$, which would play an important role in the study of higher power moments of $S(q_1q_2x; \frac{\displaystyle a_1}{\displaystyle q_1}, \frac{\displaystyle a_2}{\displaystyle q_2})$.
\begin{theorem}\label{th1}
Suppose $A\geq0$ is a fixed real number, $T\gg(q_1q_2)^\varepsilon$, then
\begin{equation*}
    \int_1^T|\Delta(q_1q_2x; r_1, q_1, r_2, q_2)|^Adx\ll\left\{\begin{array}{ll}T^{1+\frac{A}{4}}\mathcal{L}^{4A},&\text{if }~ 0\leq A\leq\frac{262}{27};\\T^{\frac{131A+154}{416}}\mathcal{L}^{4A+2},&\text{if }~  A >\frac{262}{27}. \end{array} \right.
\end{equation*}
\end{theorem}
\begin{remark}
Let $A_0>1$ be a constant such that
\begin{equation}\label{pre}
  \int_1^{T}{|\Delta(q_1q_2x; r_1, q_1, r_2, q_2)|^{A_0}}dx\ll T^{1+\frac{A_0}{4}+\varepsilon}.
\end{equation}
Then Theorem \ref{th1} shows that we can take $A_0=\frac{262}{27}$. The value $\frac{262}{27}$, which is obtained by using \eqref{9a}, is closely related to the upper bound of $\Delta(q_1q_2x; r_1, q_1, r_2, q_2)$. It is easy to prove that the conjecture \eqref{conjecture} is equivalent to \eqref{pre} holding for any $A_0>0$.
\end{remark}

 For $S\big(q_1q_2x; \frac{a_1}{q_1}, \frac{a_2}{q_2}\big)$, we have the following results.

\begin{theorem} \label{upper bound}$($\textbf{Upper bound}$)$ For $x\geq1, 1\leq a_i\leq q_i$ and $(a_i, q_i)=1$ $(i=1, 2)$, we have
\begin{equation*}
  S\big(q_1q_2x; \frac{a_1}{q_1}, \frac{a_2}{q_2}\big)\ll q_1q_2x^\frac{131}{416}(\log x)^\frac{26947}{8320}.
\end{equation*}
\end{theorem}
\begin{remark}Conjecture \eqref{conjecture} implies that
\begin{equation*}
  S\big(q_1q_2x; \frac{a_1}{q_1}, \frac{a_2}{q_2}\big)\ll q_1q_2x^{\frac{1}{4}+\varepsilon}.
\end{equation*}
\end{remark}

\begin{theorem}\label{one moment}$($\textbf{Mean value}$)$
For a large number $T$, if $1\leq a_i\leq q_i$ and $(a_i, q_i)=1$ $(i=1, 2)$, then
  \begin{equation*}
   \int_1^T{S\big(q_1q_2x; \frac{a_1}{q_1}, \frac{a_2}{q_2}\big)}dx\ll q_1q_2T^{\frac{3}{4}}.
\end{equation*}
\end{theorem}
\begin{theorem}\label{thmmain2}$($\textbf{Mean-square}$)$
Suppose $T\gg (q_1q_2)^{\varepsilon'}$ is large enough, $q_1\geq2$, $q_2\geq 3$, $1\leq a_i\leq q_i$ and $(a_i, q_i)=1$ $(i=1, 2)$. Then
\begin{equation*}
  \int_1^{T}{S^2\big(q_1q_2x; \frac{a_1}{q_1}, \frac{a_2}{q_2}\big)}dx
  =\frac{(q_1q_2)^2}{2^6\pi^2}B_2(\frac{a_1}{q_1}, \frac{a_2}{q_2}) \int_1^Tx^{\frac{1}{2}}dx+O\big((q_1q_2)^2T^{\frac{5}{4}}\mathcal{L}^3\big),
\end{equation*}
where  the constant $$B_2\big(\frac{a_1}{q_1}, \frac{a_2}{q_2}\big)=\sum_{n=1}^{\infty}\frac{\big(\Delta d_2(n;a_1, q_1, a_2, q_2)\big)^2}{n^{\frac{3}{2}}}\asymp1.$$

\end{theorem}
\medskip

\begin{theorem}\label{thmmaink}$($\textbf{Higher-power moments}$)$
Let $T\gg (q_1q_2)^{\varepsilon'}$ be a large number, $q_1\geq2$, $q_2\geq 3$, \mbox{$1\leq r_i\leq q_i$} and $(a_i, q_i)=1$ $(i=1, 2)$. If $A_{0}>9$ satisfies \eqref{pre}, then for any fixed integer $3\leq k<A_0$, we have
\begin{align}\label{asyk}
  &\int_1^{T}{S^k\big(q_1q_2x; \frac{a_1}{q_1}, \frac{a_2}{q_2}\big)}dx\\
  \nonumber=&\frac{(q_1q_2)^k}{2^{\frac{7}{2}k-1}\pi^k}B_k(\frac{a_1}{q_1}, \frac{a_2}{q_2}) \int_1^Tx^{\frac{k}{4}}dx+O\big((q_1q_2)^{k} T^{1+\frac{k}{4}-\frac{A_0-k}{4(A_0-2)s(K_0)}+\varepsilon}\big),
\end{align}
where
\begin{align*}
B_k(\frac{a_1}{q_1}, \frac{a_2}{q_2})
=&\!\sum_{\upsilon=1}^{k-1}\!\cos\!\Big(\!\frac{3\pi(k-2\upsilon)}{4}\!\Big)\!{k-1\choose \upsilon}\sum_{\{n_1,\ldots n_k\}\in\mathcal{S}(k,\upsilon)}{\prod_{i=1}^{k}\!\frac{\Delta d_2(n_i;a_1, q_1, a_2, q_2)}{n_i^\frac{3}{4}}},
  \end{align*}
with $\Delta d_2(n;a_1, q_1, a_2, q_2)$, $K_0, s(K_0)$ and $\mathcal{S}(k,\upsilon)$ defined in Section~\ref{yubei}. Clearly, $B_k\asymp1$.
\end{theorem}
\begin{remark}
From Theorem \ref{th1} we can get the asymptotic formula \eqref{asyk} for $k\in\{3,4,\ldots,9\}$, by taking $A_0=\frac{262}{27}$. Moreover, \eqref{asyk} provides the exact form of the main term for $\int_1^T{S^k\big(q_1q_2x; \frac{a_1}{q_1}, \frac{a_2}{q_2}\big)}dx$ with $k\geq10$, and if the conjecture \eqref{conjecture} were true, then we could get an asymptotic formula for $\int_1^T{S^k\big(q_1q_2x; \frac{a_1}{q_1}, \frac{a_2}{q_2}\big)}dx$ with any $k\geq3$.
 \end{remark}
 If $A_0=\frac{262}{27}$, then $K_0=10$, and $s(K_0)=257$. Hence, we get an immediate corollary of Theorem \ref{thmmaink}.
\begin{cor}
Let $3\leq k\leq9$ be a fixed integer. If $T, a_i$ and $q_i (i=1,2)$ satisfying the hypothesis of Theorem \ref{thmmaink}, then
\begin{align*}
  &\int_1^{T}\!\!\!\!{S^k\big(q_1q_2x; \frac{a_1}{q_1}, \frac{a_2}{q_2}\big)}dx\!=\!\frac{(q_1q_2)^k}{2^{\frac{7}{2}k-\!1}\pi^k}B_k(\frac{a_1}{q_1}, \frac{a_2}{q_2}) \int_1^T\!\!\!\!x^{\frac{k}{4}}dx\!+\!O\big((q_1q_2)^{k} T^{1+\frac{k}{4}\!-\!\frac{262-27k}{213824}\!+\varepsilon}\big).
\end{align*}
\end{cor}
\begin{remark} There are similar conclusions for $\mathop{{\sum}'}_{mn\leq x}{\sin(2\pi n\theta_1)\sin(2\pi m\theta_2)}$ and $\mathop{{\sum}'}_{mn\leq x}{\cos(2\pi n\theta_1)\cos(2\pi m\theta_2)}$, which we do not present in this paper.
\end{remark}
\textbf{Acknowledgement.}The authors deeply thank the referee for carefully reading the manuscript and valuable suggestions.
\vspace{1ex}
\section{Background of divisor problem with congruence conditions}

In this section, we review some background of the divisor problem with congruence conditions.

 \subsection{Some estimates for $\Delta(x; r_1, q_1, r_2, q_2)$}

 It follows from Huxley's estimates\cite{Huxley03} that
\begin{equation}\label{9a}
  \Delta(x; r_1, q_1, r_2, q_2)\ll\Big(\frac{x}{q_1q_2}\Big)^\frac{131}{416}\Big(\log\Big(\frac{x}{q_1q_2}\Big)\Big)^\frac{26947}{8320}
\end{equation}
uniformly in $1\leq r_1\leq q_1\leq x, 1\leq r_2\leq q_2\leq x$. It is conjectured that
\begin{equation}\label{conjecture}
\Delta(x; r_1, q_1, r_2, q_2)\ll\Big(\frac{x}{q_1q_2}\Big)^{\frac{1}{4}+\varepsilon}
\end{equation}
uniformly in $1\leq r_1\leq q_1\leq x, 1\leq r_2\leq q_2\leq x$, $\forall\varepsilon>0$, which is an analogue of the well-known conjecture that $\Delta(x)\ll x^{\frac{1}{4}+\varepsilon}$.

M\"{u}ller and Nowak\cite{MullerNowak} studied the mean value of $\Delta(x; r_1, q_1, r_2, q_2)$. They pointed out
\begin{equation}\label{muller}
  \int_1^T{\Delta(x; r_1, q_1, r_2, q_2)}dx\ll (q_1q_2)^{-\frac{1}{4}}T^{\frac{3}{4}}
\end{equation}
and
\begin{equation}\label{muller2}
 \int_1^T{\Delta^2(x; r_1, q_1, r_2, q_2)}dx=c_2(q_1q_2)^{\frac{1}{2}}T^\frac{3}{2}+o\big((q_1q_2)^{\frac{1}{2}}T^\frac{3}{2}\big)
\end{equation}
uniformly in $1\leq r_i\leq q_i\leq T$ $(i=1, 2)$, if $T$ is a large number, and $c_2$ is a constant.

\subsection{An analogue Voronoi formula for $\Delta(x; r_1, q_1, r_2, q_2)$}

In this subsection we shall prove the following.
\begin{lemma}\label{lem:vor}
Let $J=[\frac{\mathcal{L}+2\log q_1q_2-4\log \mathcal{L}}{\log2}]$, $H\geq2$ be a parameter to be determined, and $ T^\varepsilon<y\leq \min(H^2, (q_1q_2)^2 T)\mathcal{L}^{-4}$. Suppose $\frac{T}{2}\leq x\leq T$. Denote
\begin{align*}
  &\tau(n, x)=\sum_{n=hl}\cos\bigg(4\pi \sqrt{nx}-2\pi\Big(\frac{hr_2}{q_2}+\frac{lr_1}{q_1}+\frac{1}{8}\Big)\bigg), \\
  &\tau_{12}(n, x; H) =\mathop{{\sum}'}_
{\{h,l\}\in E(n; H,J)}
\cos\bigg(4\pi \sqrt{nx}-2\pi\Big(\frac{hr_2}{q_2}+\frac{lr_1}{q_1}+\frac{1}{8}\Big)\bigg),\\
  & \tau_{21}(n, x; H) =\mathop{{\sum}'}_
{\{h,l\}\in E(n; H,J)}
\cos\bigg(4\pi \sqrt{nx}-2\pi\Big(\frac{hr_1}{q_1}+\frac{lr_2}{q_2}+\frac{1}{8}\Big)\bigg),
\end{align*}
with $E(n; H,J)$ defined in Section~\ref{yubei}. Then
\begin{multline}\label{lemliu}
  \Delta(q_1q_2x; r_1, q_1, r_2, q_2)=R'_0(x; y)+{R_{12}}'(x; y, H)+{R_{21}}'(x; y, H)\\+G'_{12}(x; H)+G'_{21}(x; H)
  +O\big(\log^3(q_1q_2T)\big).
\end{multline}
Here
\begin{eqnarray*}
  \nonumber &&R'_0(x; y)=\frac {x^{\frac{1}{4}}}{\sqrt{2}\pi}\sum_{n\leq y}\frac{\tau(n, x)}{n^{\frac{3}{4}}},  \\
  \nonumber &&{R_{12}}'(x; y, H)=\frac {x^{\frac{1}{4}}}{\sqrt{2}\pi}\sum_{y<n\leq 2^{J+1}H^2}\frac{\tau_{12}(n, x; H)}{n^{\frac{3}{4}}},  \\
  &&{R_{21}}'(x; y, H)=\frac {x^{\frac{1}{4}}}{\sqrt{2}\pi}\sum_{y<n\leq 2^{J+1}H^2}\frac{\tau_{21}(n, x; H)}{n^{\frac{3}{4}}},  \\
 \nonumber &&G'_{12}(x; H)=\sum_{\begin{subarray}{c}  n_1\leq q_1\sqrt{T}\\ n_1\equiv r_1\!\!\!\!\!\pmod{q_1}\end{subarray}}O\biggl(\min\Bigl(1, \frac{1}{H\|\frac{q_1x}{n_1}-\frac{r_2}{q_2}\|}\Bigr)\biggr),\\
  &&G'_{21}(x; H)=\sum_{\begin{subarray}{c}  n_2\leq q_2\sqrt{T}\\ n_2\equiv r_2\!\!\!\!\!\pmod{q_2}\end{subarray}}O\bigg(\min\Big(1, \frac{1}{H\|\frac{q_2x}{n_2}-\frac{r_1}{q_1}\|}\Big)\bigg).
\end{eqnarray*}
\end{lemma}

\paragraph{Proof of Lemma \ref{lem:vor}}

From \eqref{4}, we write
\begin{equation}\label{00}
    \Delta(q_1q_2x; r_1, q_1, r_2, q_2)=F_{1}+F_{2}+O(1),
\end{equation}
where
\begin{align*}
    &F_{1}(x)=-\!\!\!\!\sum_{\begin{subarray}{c}n_1\leq q_1\sqrt{x}\\n_1\equiv r_1\!\!\!\!\!\pmod{\!q_1}\end{subarray}}\!\psi\big(\frac{q_1x}{n_1}-\frac{r_2}{q_2}\big),
    & F_{2}(x)=-\!\!\!\!\sum_{\begin{subarray}{c}n_2\leq q_2\sqrt{x}\\n_2\equiv r_2\!\!\!\!\!\pmod{\!q_2}\end{subarray}}\!\psi\big(\frac{q_2x}{n_2}-\frac{r_1}{q_1}\big).
\end{align*}

The rest of the proof is almost the same as that for Lemma 3.1 in Liu \cite{LiuKui}.\qed

\vspace{2ex}

\section{Preliminaries}\label{yubei}

We introduce some notations first.
Set
\begin{align*}
    E(n; H,J)=\{\{h,l\}\mid hl=n,1\leq h\leq H,h\leq l\leq2^{J+1}h\}.
\end{align*}
Denote
\begin{align*}
\Delta d_2(n;a_1, q_1, a_2, q_2)
=&d(n;  a_1, q_1, a_2, q_2)+d(n; - a_1, q_1, a_2, q_2)\\
&-d(n;  a_1, q_1, -a_2, q_2)-d(n; - a_1, q_1, a_2, q_2),
\end{align*}
\begin{align*}
&\Delta d_{2,1}(n,H,J;a_1, q_1, a_2, q_2)\\
=&\!\!\!\!\mathop{{\sum}'}_{\begin{subarray}{c}\{h,l\}\in E(n; H,J)\\ h\equiv a_2\!\!\!\!\!\pmod{q_2}\\l\equiv  a_1\!\!\!\!\!\pmod{q_1}\end{subarray}}\!\!\!1+\!\!\!\!\mathop{{\sum}'}_{\begin{subarray}{c}\{h,l\}\in E(n; H,J)\\ h\equiv a_2\!\!\!\!\!\pmod{q_2}\\l\equiv - a_1\!\!\!\!\!\pmod{q_1}\end{subarray}}\!\!\!1-\!\!\!\mathop{{\sum}'}_{\begin{subarray}{c}\{h,l\}\in E(n; H,J)\\ h\equiv -a_2\!\!\!\!\!\pmod{q_2}\\l\equiv a_1\!\!\!\!\!\pmod{q_1}\end{subarray}}\!\!\!1-\!\!\!\mathop{{\sum}'}_{\begin{subarray}{c}\{h,l\}\in E(n; H,J)\\ h\equiv -a_2\!\!\!\!\!\pmod{q_2}\\l\equiv - a_1\!\!\!\!\!\pmod{q_1}\end{subarray}}\!\!\!1,
\end{align*}
\begin{align*}
&\Delta d_{2,2}(n,H,J;a_1, q_1, a_2, q_2)\\
=&\!\!\!\!\mathop{{\sum}'}_{\begin{subarray}{c}\{h,l\}\in E(n; H,J)\\h\equiv  a_1\!\!\!\!\!\pmod{q_1}\\ l\equiv a_2\!\!\!\!\!\pmod{q_2}\end{subarray}}\!\!\!1+\!\!\!\!\mathop{{\sum}'}_{\begin{subarray}{c}\{h,l\}\in E(n; H,J)\\h\equiv - a_1\!\!\!\!\!\pmod{q_1}\\ l\equiv a_2\!\!\!\!\!\pmod{q_2}\end{subarray}}\!\!\!1-\!\!\!\mathop{{\sum}'}_{\begin{subarray}{c}\{h,l\}\in E(n; H,J)\\h\equiv a_1\!\!\!\!\!\pmod{q_1}\\ l\equiv -a_2\!\!\!\!\!\pmod{q_2}\end{subarray}}\!\!\!1-\!\!\!\mathop{{\sum}'}_{\begin{subarray}{c}\{h,l\}\in E(n; H,J)\\h\equiv - a_1\!\!\!\!\!\pmod{q_1}\\ l\equiv -a_2\!\!\!\!\!\pmod{q_2}\end{subarray}}\!\!\!1.
\end{align*}
Set $\textbf{1}=\{0, 1\}$, $\textbf{\textit{i}}=(i_1, \cdots, i_{k-1})\in \textbf{1}^{k-1}$, $\textbf{\textit{n}}=(n_1, \cdots, n_{k})\in \mathbb{N}^{k}$. Define
\begin{align*}
   & |\textbf{\textit{i}}|=i_1+\cdots+i_{k-1},\\
    &\alpha(\textbf{\textit{n}}; \textbf{\textit{i}})=\sqrt{n_1}+(-1)^{i_1}\sqrt{n_2}+\cdots+(-1)^{i_{k-1}}\sqrt{n_k}, \\
   & \beta(\textbf{\textit{i}})=
    \frac{3\pi}{4}+(-1)^{i_1}\frac{3\pi}{4}+\cdots+(-1)^{i_{k-1}}\frac{3\pi}{4}=(k-2|\textbf{\textit{i}}|)\frac{3\pi}{4},\\
    &\mathcal{S}(k,\upsilon)=\{\{n_1,\ldots n_k\}\mid\sqrt{n_1}\!+\!\cdots\!+\!\sqrt{n_\upsilon}\!=\!\sqrt{n_{\upsilon+1}}\!+\!\cdots\!+\!\sqrt{n_k}\}.
    \end{align*}
Suppose $A_0>2$. Denote
\begin{align*}
  &K_0:=\min\{n\in\mathbb{N}:n\geq A_0, 2|n\},
  &&s(k):=2^{k-2}+\frac{k-6}{4}.
\end{align*}

To prove the theorems, we need the following definition and Lemmas.
\begin{definition}\label{expa}$($ See eg.\cite[p.73]{ivic} $)$
Let $A>\frac{1}{2}$, $B\geq 1$, and $f(x)\in C^5[B,2B]$ be real valued function. Assume the derivatives of $f(x)$ for $B\leq x\leq 2B$  satisfy
$$AB^{1-r}\ll|f^{(r)}(x)|\ll AB^{1-r},\quad r=1,2,3,4,$$
where the $\ll$ constants depend on $r$ alone.  Consider
\begin{equation*}
    S_f(d)=\sum_{B<n\leq B+d}e(f(n)),\quad 1<d\leq B.
\end{equation*}
The pair of nonnegative real numbers $(\kappa, \lambda)$ will be called an exponent pair,
if $0\leq\kappa \leq\frac{1}{2}\leq\lambda\leq 1$ and
$$S_f(d)\ll A^\kappa B^\lambda$$
holds uniformly for $f$ and $d$.
\end{definition}
\begin{lemma}\label{lem:Vaa}$($ See \cite[p. 210]{vaaler1985some} $)$
Let $H\geq2$ be any real number. Then
\begin{equation*}
\psi(u)=\sum_{1\leq |h|\leq H}a(h)e(hu)+O\big(\sum_{1\leq |h|\leq H}b(h)e(hu)\big),
\end{equation*}
where $a(h)$ and $b(h)$ are functions such that $a(h)\ll\frac{1}{|h|}$ and $b(h)\ll\frac{1}{H}$.
\end{lemma}
\begin{lemma}\label{lem:min}$($ See \cite[Theorem 2.2.2)]{Min1981method} $)$
Let $A_1,\ldots,A_5$ be absolute positive constants, $f(x)\in C^3[a,b]$ and $g(x)\in C^1[a,b]$  be $\mathbb{C}$ valued functions satisfying
\begin{align*}
    &\frac{A_1}{R}\leq|f''(x)|\leq\frac{A_1}{R}, &&|f'''(x)|\leq\frac{A_3}{RU},\ U\geq 1,\\
    &|g(x)|\leq A_4G,&&|g'(x)|\leq \frac{A_5G}{U_1},\ \ U_1\geq 1.
\end{align*}
Suppose $\alpha\leq f'(x)\leq \beta$ for $a \leq x\leq b$. Then
\begin{align*}
    \sum_{a < x\leq b}g(n)e\big(f(n)\big)=&\sum_{\alpha\leq u\leq \beta}b_u\frac{g(n_u)}{\sqrt{f''(n_u)}}e\big(f(n_u)-un_u+\frac{1}{8}\big)\\
    &+O\big(G\log(\beta-\alpha+2)+G(b-a+R)(U^{-1}+U_1^{-1})\big)\\
    &+O\Big(G\min\Big\{\sqrt{R}, \max\Big(\frac{1}{\langle\alpha\rangle},\frac{1}{\langle\beta\rangle}\Big)\Big\}\Big),
\end{align*}
where $n_u$ is the solution of $f'(n)=u$,
\begin{align*}
   &\langle t\rangle =\left\{\begin{array}{ll}
   \|t\|,&\text{if }t\text{  is not an integer},\\
   \beta-\alpha,& \text{if }t\text{  is an integer},
   \end{array}\right.\\
   &b_u=\left\{\begin{array}{ll}
   1,&\text{if }\alpha< u<\beta,\text{ or } \alpha,~\beta \text{ are not integers},\\
   \frac{1}{2},& \text{if }u=\alpha \text{ or}~\beta\text{  is an integer},
   \end{array}\right.\\
   &\sqrt{f''}=\left\{\begin{array}{ll}
   \sqrt{f''},&\text{if }f''>0,\\
   i\sqrt{f''},& \text{if }f''<0.
   \end{array}\right.\\
\end{align*}
\end{lemma}
\begin{lemma}\label{lem:ivic}$($Hal\'{a}sz-Montgomery inequality$)$ $($see e.g., Ivi\'{c}\cite{ivic}$)$ Let $\xi$, $\varphi_1$,\ldots, $\varphi_R$ are arbitrary vectors in an inner-product vector space over $\mathbb{C}$, where $(a,b)$ will be the notion for the inner product and $\| a\|^2=(a,a)$. Then
\begin{equation*}
    \sum_{r\leq R}\big|(\xi,\varphi_r)\big|^2\leq \| \xi\|^2\max_{r\leq R}\sum_{s\leq R}\big|(\varphi_r,\varphi_s)\big|.
\end{equation*}
\end{lemma}
\begin{lemma}\label{lemsum}$($See \cite[eq.(4.6)]{zhai04} $)${\, \, \rm} Suppose $k\geq2$ is an integer, $y$ is a large positive number, $\textbf{\textit{i}}=(i_1, \cdots, i_{k-1})\in \textup{\textbf{1}}^{k-1}$, and $\textbf{\textit{n}}=(n_1, \cdots, n_{k})\in \mathbb{N}^{k}$. Then
\begin{equation*}
  \sum_{\textbf{\textit{i}}\in \textup{\textbf{1}}^{k-1}}\sum_{\begin{subarray}{c}n_i\leq y\\1\leq i\leq k\\\alpha(\textbf{\textit{n}}; \textbf{\textit{i}})\neq0\end{subarray}}\frac{d(n_1)\cdots d(n_k)}{(n_1\cdots n_k)^\frac{3}{4} \alpha(\textbf{\textit{n}};\textbf{\textit{i}})}\ll y^{s(k)+\varepsilon}.
\end{equation*}
\end{lemma}

\begin{lemma}\label{lem2.3}$($ See \cite[Lemma 3.1]{zhai04} $)${\, \, \rm}Suppose $f(n):\mathbb{N}\rightarrow\mathbb{R}$ is a function with $f(n)\ll n^{\varepsilon^\prime}$, $k\geq2$ is an integer, $1\leq\upsilon<k$ is a fixed integer, $y>1$ is large. Define
\begin{align*}
  s_{k, \upsilon}(f; y):&=\sum_{\begin{subarray}{c}\{n_1,n_2,\cdots, n_k\}\in\mathcal{S}(k,\upsilon)\\n_1,n_2,\cdots, n_k\leq y\end{subarray}}\frac{f(n_1)f(n_2)\cdots f(n_k)}{(n_1n_2\cdots n_k)^{\frac{3}{4}}},\quad 1\leq\upsilon<k,\\
  s_{k, \upsilon}(f):&=\sum_{\{n_1,n_2,\cdots, n_k\}\in\mathcal{S}(k,\upsilon)}\frac{f(n_1)f(n_2)\cdots f(n_k)}{(n_1n_2\cdots n_k)^{\frac{3}{4}}},\quad  1\leq\upsilon<k.
\end{align*}
  Then  for any $\varepsilon>0$, $1\leq\upsilon<k$, we have
    \begin{gather*}
    s_{k, \upsilon}(f)\ll 1,\\
   |s_{k, \upsilon}(f)-s_{k, \upsilon}(f; y)|\ll y^{-\frac{1}{2}+\varepsilon}.
  \end{gather*}
\end{lemma}

\begin{lemma}\label{lem2.4}$($ Hilbert's inequality $)$$($ See e.g.\cite{shan} $)${\, \, \rm} Let $x_1<x_2<\cdots<x_n$ be a sequence of real numbers.  If there exists$ ~\delta>0$, such that $\min\limits_{s\neq r}|x_r-x_s|\geq \delta_r\geq\delta>0 (1\leq r\leq n)$, then there exists an absolute constant $C$, such that
\begin{equation*}
  \bigg|\sum_{s\neq r}u_r\bar{{u_s}}(x_r-x_s)^{-1}\bigg|\leq C\sum_r {\delta_r}^{-1}\left|u_r\right|^2,
\end{equation*}
for arbitrary complex numbers $u_1,u_2,\cdots,u_n$.
\end{lemma}

\vspace{1ex}

\section{ Proofs of Theorem \ref{th1}}

We follow the approach of Zhai \cite{zhai2008error}, and first present a large value estimate of $\Delta(q_1q_2x; r_1, q_1, r_2, q_2)$.

\subsection{Large value estimate of $\Delta(q_1q_2x; r_1, q_1, r_2, q_2)$}

We shall prove the following
\begin{theorem}\label{th:large}
Suppose $T\gg (q_1q_2)^\varepsilon$. Let $\frac{T}{2}\leq x_1<\cdots<x_M\leq T$ satisfy $|\Delta(q_1q_2x_s; r_1, q_1, r_2, q_2)|\gg V$ $(s=1,\ldots,M)$ and $|x_i-x_j|\geq V\gg T^\frac{7}{32}\mathcal{L}^4 (i\neq j)$. Then
\begin{equation*}
    M\ll TV^{-3}\mathcal{L}^9+T^\frac{15}{4}V^{-12}\mathcal{L}^{41}.
\end{equation*}
\end{theorem}

\paragraph{Proof of Theorem \ref{th:large}}
Let $T_0>V$ be a parameter to be determined and $I$ be any subinterval of $[\frac{T}{2},T]$  of length not exceeding $T_0$. Without loss of generality, we may assume $G=I\cap\{x_1,\ldots,x_M\}=\{x_1,\ldots,x_{M_0}\}$.

From \eqref{00}, let $J_0=\big[\frac{\mathcal{L}+2\log\mathcal{L}-2\log V}{\log 2}\big]$ and
\begin{align*}
&\Psi_1(x;j)=\sum_{\begin{subarray}{c}k_1\\2^{-\frac{j+1}{2}}q_1\sqrt{x}\leq k_1 q_1+r_1\leq 2^{-\frac{j}{2}}q_1\sqrt{x}\end{subarray}}\!\!\psi\Big(\frac{q_1x}{k_1 q_1+r_1}-\frac{r_2}{q_2}\Big),\\
&\Psi_2(x;j)=\sum_{\begin{subarray}{c}k_2\\2^{-\frac{j+1}{2}}q_2\sqrt{x}\leq k_2 q_2+r_2\leq 2^{-\frac{j}{2}}q_2\sqrt{x}\end{subarray}}\!\!\psi\Big(\frac{q_2x}{k_2 q_2+r_2}-\frac{r_1}{q_1}\Big).
\end{align*}
Then
\begin{equation}\label{psi0}
  \Delta(q_1q_2x; r_1, q_1, r_2, q_2)=-\sum_{j=0}^{J_0}\Psi_1(x;j)-\sum_{j=0}^{J_0}\Psi_2(x;j)+O\Big(\frac{V}{\mathcal{L}}\Big).
\end{equation}
By Cauchy's inequality we get
\begin{equation*}
    \Delta^2(q_1q_2x; r_1, q_1, r_2, q_2)\ll\mathcal{L}\sum_{j=0}^{J_0}\Psi_1^2(x;j) +\mathcal{L}\sum_{j=0}^{J_0}\Psi_2^2(x;j)+O\Big(\frac{V^2}{\mathcal{L}^2}\Big),
\end{equation*}
which implies that for $x\in G$
\begin{equation*}
    \Delta^2(q_1q_2x; r_1, q_1, r_2, q_2)\ll\mathcal{L}\sum_{j=0}^{J_0}\Psi_1^2(x;j) +\mathcal{L}\sum_{j=0}^{J_0}\Psi_2^2(x;j).
\end{equation*}

Denote
\begin{align*}
    &G_1:=\Big\{x_s\in G~\Big|~\sum_{j=0}^{J_0}\Psi_1^2(x_s;j)\geq\sum_{j=0}^{J_0}\Psi_2^2(x_s;j)\Big\},&M_1=:\#G_1;\\
    &G_2:=\Big\{x_s\in G~\Big|~\sum_{j=0}^{J_0}\Psi_2^2(x_s;j)\geq\sum_{j=0}^{J_0}\Psi_1^2(x_s;j)\Big\},&M_2=:\#G_2.
\end{align*}

We first estimate $M_1$. Assume that $G_1=\{x_1,\ldots,x_{M_1}\}$. Then we have
\begin{equation*}
    \Delta^2(q_1q_2x_s; r_1, q_1, r_2, q_2)\ll\mathcal{L}\sum_{j=0}^{J_0}\Psi_1^2(x_s;j),\quad \text{for~} 1\leq s\leq M_1.
\end{equation*}
Summing  over the set $G_1$ we get
\begin{align}\label{04}
    M_1V^2 \ll& \sum_{s\leq M_1}\Delta^2(q_1q_2x_s; r_1, q_1, r_2, q_2)\\
   \nonumber  \ll&\mathcal{L}\sum_{j=0}^{J_0}\sum_{s\leq M_1}\Psi_1^2(x_s;j) \ll \mathcal{L}^2\sum_{s\leq M_1}\Psi_1^2(x_s;j),
\end{align}
for some fixed $j$ with $0\leq j\leq J_0$. Let $2\leq H<T$ be a parameter to be determined later. From Lemma \ref{lem:Vaa} and Lemma \ref{lem:min}, we have
\begin{align*}
    \Psi_1(x_s;j)\ll \bigg|\sum_{1\leq h\leq H}c(h)\!\!\!\!\sum_{\begin{subarray}{c}k_1\\2^{-\frac{j+1}{2}}q_1\sqrt{x_s}\leq k_1 q_1+r_1\leq 2^{-\frac{j}{2}}q_1\sqrt{x_s}\end{subarray}}\!\!\!\!e\Big(-h\big(\frac{q_1x_s}{k_1 q_1+r_1}-\frac{r_2}{q_2}\big)\Big)\bigg|,
\end{align*}
where $c(h)$ is some function such that $c(h)\ll \frac{1}{h}$.
By using Lemma \ref{lem:min}, we get
\begin{align}\label{05}
\Psi_1(x_s;j)\!\ll\!\frac{x^\frac{1}{4}}{\sqrt{2}}\! \sum_{1\leq h\leq H}\!\!c(h)\!\!\!\sum_{2^{j}h\leq l\leq2^{j+1}h}\!\!\frac{b_lh^\frac{1}{4}}{l^\frac{3}{4}}e\big(\!-\!2 \sqrt{hlx_s}\!+\!\frac{hr_2}{q_2}\!+\!\frac{lr_1}{q_1}\big)\!+\!O(\mathcal{L}^2).
\end{align}
Inserting \eqref{05} into \eqref{04} yields
\begin{align*}
    M_1V^2 \!\ll&\mathcal{L}^2T^\frac{1}{2}\sum_{s\leq M_1}\Big|\sum_{1\leq h\leq H}\!\!c(h)\!\!\!\sum_{2^{j}h\leq l\leq2^{j+1}h}\!\!\frac{b_lh^\frac{1}{4}}{l^\frac{3}{4}}e\big(\!-\!2 \sqrt{hlx_s}\!+\!\frac{hr_2}{q_2}\!+\!\frac{lr_1}{q_1}\big)\Big|^2\!\!+\!\!M_1\mathcal{L}^6\\
    \ll&\mathcal{L}^2T^\frac{1}{2}\sum_{s\leq M_1}\Big|\sum_{1\leq n\leq2^{j+1}H^2}\sum_{\begin{subarray}{c}n=hl\\1\leq h\leq H\\2^{j}h\leq l\leq2^{j+1}h\end{subarray}}\frac{b_lh c(h)}{n^\frac{3}{4}}e\big(\!-\!2 \sqrt{hlx_s}+\frac{hr_2}{q_2}+\frac{lr_1}{q_1}\big)\Big|^2.
\end{align*}
Let
\begin{align*}
    \gamma(n;H,j)=\sum_{\begin{subarray}{c}n=hl\\1\leq h\leq H\\2^{j}h\leq l\leq2^{j+1}h\end{subarray}}b_lh c(h)e\big(\frac{hr_2}{q_2}+\frac{lr_1}{q_1}\big).
\end{align*}
Then  it is easy to see that $\gamma(n;H,j)\ll d(n)$, and
\begin{align*}
    M_1V^2 \ll&\mathcal{L}^2T^\frac{1}{2}\sum_{s\leq M_1}\Big|\sum_{1\leq n\leq2^{j+1}H^2}\frac{\gamma(n;H,j)}{n^\frac{3}{4}}e(-2 \sqrt{nx_s})\Big|^2.
\end{align*}
By a splitting argument and  Cauchy's inequality, we get
\begin{align*}
    M_1V^2\ll&\mathcal{L}^2T^\frac{1}{2}\sum_{s\leq M_1}\Big|\sum_\nu\sum_{2^{j-v}H^2< n\leq2^{j+1-v}H^2}\frac{\gamma(n;H,j)}{n^\frac{3}{4}}e(-2 \sqrt{nx_s})\Big|^2\\
 \nonumber   \ll&\mathcal{L}^3T^\frac{1}{2}\sum_\nu\sum_{s\leq M_1}\Big|\sum_{2^{j-v}H^2< n\leq2^{j+1-v}H^2}\frac{\gamma(n;H,j)}{n^\frac{3}{4}}e(-2 \sqrt{nx_s})\Big|^2\\
  \nonumber  \ll&\mathcal{L}^4T^\frac{1}{2}\sum_{s\leq M_1}\Big|\sum_{2^{j-v}H^2< n\leq2^{j+1-v}H^2}\frac{\gamma(n;H,j)}{n^\frac{3}{4}}e(-2 \sqrt{nx_s})\Big|^2
\end{align*}
for some $0\leq \nu\ll\mathcal{L}$. Let $N_1=2^{j-v}H^2$. Take $\xi=\{\xi_n\}_{n=1}^\infty$ with $\xi_n=\gamma(n;H,j)n^{-\frac{3}{4}}$ for $N_1< n\leq2N_1$ and zero otherwise, and $\varphi_s=\{\varphi_{s,n}\}_{n=1}^\infty$ with $\varphi_{s,n}=e(-2 \sqrt{nx_s})$ for $N_1< n\leq2N_1$ and zero otherwise. Then
\begin{gather*}
    \sum_{s\leq M_1}(\xi,\varphi_{s})^2=\sum_{s\leq M_1}\Big|\sum_{N_1< n\leq2N_1}\frac{\gamma(n;H,j)}{n^\frac{3}{4}}e(-2 \sqrt{nx_s})\Big|^2,\\
    (\varphi_{s_1},\varphi_{s_2})=\sum_{N_1< n\leq2N_1}e\big(2 \sqrt{n}(\sqrt{x_{s_2}}-\sqrt{x_{s_1}})\big),\\
    \|\xi\parallel^2=\sum_{N_1< n\leq2N_1}\frac{|\gamma(n;H,j)|^2}{n^\frac{3}{2}}\ll\sum_{N_1< n\leq2N_1}\frac{d^2(n)}{n^\frac{3}{2}}\ll N_1^{-\frac{1}{2}}\mathcal{L}^3,
\end{gather*}
where we used the estimate $\sum_{n\leq N}d^2(n)\ll N\log^3N$. Thus by Lemma \ref{lem:ivic}, we get
\begin{align}\label{07}
    M_1V^2\!\ll&T^\frac{1}{2}N_1^{-\frac{1}{2}}\mathcal{L}^7\max_{s_1\leq M_1}\sum_{s_2\leq M_1}\Big|\sum_{N_1< n\leq2N_1}e\big(2 \sqrt{n}(\sqrt{x_{s_2}}-\sqrt{x_{s_1}})\big)\Big|\\
 \nonumber= &T^\frac{1}{2}N_1^{\frac{1}{2}}\mathcal{L}^7\!+\!T^\frac{1}{2}N_1^{\!-\!\frac{1}{2}}\mathcal{L}^7\! \max_{s_1\leq M_1}\sum_{\begin{subarray}{c}s_2\leq M_1\\s_2\neq s_1\end{subarray}}\Big|\!\sum_{N_1< n\leq2N_1}\!\!\!e\big(2 \sqrt{n}(\sqrt{x_{s_2}}\!-\!\sqrt{x_{s_1}})\big)\Big|.
\end{align}
By the Kuz'min-Landau inequality and the exponent pair $(\frac{4}{18}, \frac{11}{18})$ (see eg. \cite[p.77]{ivic}), we have
\begin{align*}
    \sum_{N_1< n\leq2N_1}e\big(2 \sqrt{n}(\sqrt{x_{s_2}}-\sqrt{x_{s_1}})\big)\ll& \frac{\sqrt{N_1}}{|\sqrt{x_{s_2}}-\sqrt{x_{s_1}}|} +\Big(\frac{|\sqrt{x_{s_2}}-\sqrt{x_{s_1}}|}{\sqrt{N_1}}\Big)^\frac{4}{18}N_1^\frac{11}{18}\\
    \ll&\frac{\sqrt{N_1T}}{|x_{s_2}-x_{s_1}|}  +\Big(\frac{|x_{s_2}-x_{s_1}|}{\sqrt{N_1T}}\Big)^\frac{4}{18}N_1^\frac{11}{18}\\
    \ll & \frac{\sqrt{N_1T}}{|s_2-s_1|V}+T^{-\frac{1}{9}}T_0^{\frac{2}{9}}N_1^\frac{1}{2},
\end{align*}
in view of $|x_{s_2}-x_{s_1}|\leq T_0$. Combining this with \eqref{07}, we get
\begin{align}\label{08}
    M_1V^2\ll&T^\frac{1}{2}N_1^{\frac{1}{2}}\mathcal{L}^7+TV^{-1}\mathcal{L}^7 \max_{s_1\leq M_1}\sum_{\begin{subarray}{c}s_2\leq M_1\\s_2\neq s_1\end{subarray}}\frac{1}{|s_2-s_1|}+M_0T^{\frac{7}{18}}T_0^{\frac{2}{9}}\mathcal{L}^7\\
   \nonumber \ll&T^\frac{1}{2}N_1^{\frac{1}{2}}\mathcal{L}^7+TV^{-1}\mathcal{L}^8 +M_0T^{\frac{7}{18}}T_0^{\frac{2}{9}}\mathcal{L}^7.
\end{align}
Take $H=T^\frac{1}{2}V^{-1}2^{-\frac{j}{2}}\mathcal{L}^2$, $T_0=T^{-\frac{7}{4}}V^{9}\mathcal{L}^{-32}$. Then $N_1=TV^{-2}2^{-\nu}\mathcal{L}^4$, and $T_0\gg V$ if $V\gg T^\frac{7}{32}\mathcal{L}^4$. Thus from \eqref{08}, we get
\begin{equation*}
    M_1\ll TV^{-3}\mathcal{L}^9.
\end{equation*}

For $M_2$, by symmetry, we can get the same estimate. Noting that $M_1+M_2\geq M_0$, we have
\begin{equation*}
    M_0\ll TV^{-3}\mathcal{L}^9.
\end{equation*}

Dividing the interval $[\frac{T}{2},T]$ into $O\big(1+\frac{T}{2T_0}\big)$ subintervals of length not exceeding $T_0=T^{-\frac{7}{4}}V^{9}\mathcal{L}^{-32}$, we see
\begin{align*}
    M\ll M_0\big(1+\frac{T}{2T_0}\big)\ll TV^{-3}\mathcal{L}^9+ T^\frac{15}{4}V^{-12}\mathcal{L}^{41}.
\end{align*}
This completes the proof of Theorem \ref{th:large}.
\qed
\subsection{Proof of Theorem \ref{th1}}
For $A=0$, this is trivial. For $0<A\leq 2$, it follows from \eqref{muller2} and H\"{o}lder's inequality. So it suffices to prove it for $A>2$.

With a similar argument to (13.70) of Ivi\'{c}, we see
\begin{equation*}
        \int_\frac{T}{2}^T|\Delta(q_1q_2x; r_1, q_1, r_2, q_2)|^Adx\ll T^{1+\frac{A}{4}}\mathcal{L}^{4A}\!+\!\sum_VV\!\sum_{s\leq N_V}|\Delta(q_1q_2x_s; r_1, q_1, r_2, q_2)|^A,
\end{equation*}
where $T^\frac{1}{4}\mathcal{L}^4\leq V\leq T^\frac{131}{416}\mathcal{L}^4$, $V<|\Delta(q_1q_2x_s; r_1, q_1, r_2, q_2)|\leq 2V$, $s=1,\ldots,N_V$, $|x_{s_1}-x_{s_2}|\geq V$ for $s_1\neq s_2\leq N_V$, and the number of $V$ is $\ll \mathcal{L}$.
By Theorem \ref{th:large}, we get
\begin{align*}
    V\sum_{s\leq N_V}|\Delta(q_1q_2x_s; r_1, q_1, r_2, q_2)|^A\ll N_VV^{A+1}\ll TV^{A-2}\mathcal{L}^9+T^\frac{15}{4}V^{A-11}\mathcal{L}^{41}.
\end{align*}
Thus we have
\begin{align*}
    \int_\frac{T}{2}^T|\Delta(q_1q_2x; r_1, q_1, r_2, q_2)|^Adx\ll& T^{1+\frac{A}{4}}\mathcal{L}^{4A}+TV^{A-2}\mathcal{L}^{10}+T^\frac{15}{4}V^{A-11}\mathcal{L}^{42},
\end{align*}
for some fixed $T^\frac{1}{4}\mathcal{L}^4\leq V\leq T^\frac{131}{416}\mathcal{L}^4$,
which implies
\begin{align*}
&\int_\frac{T}{2}^T|\Delta(q_1q_2x; r_1, q_1, r_2, q_2)|^Adx\\
\ll&\left\{\begin{array}{ll}
T^{1+\frac{A}{4}}\mathcal{L}^{4A}+T^{1+\frac{131}{416}(A-2)}\mathcal{L}^{4A+2}+T^{\frac{15}{4}+\frac{1}{4}(A-11)}\mathcal{L}^{4A-2}, &\text{if }~ 2< A\leq11,\\
T^{1+\frac{A}{4}}\mathcal{L}^{4A}+T^{1+\frac{131}{416}(A-2)}\mathcal{L}^{4A+2}+T^{\frac{15}{4}+\frac{131}{416}(A-11)}\mathcal{L}^{4A-2}, &\text{if }~  A>11, \end{array} \right.\\
    \ll&\left\{\begin{array}{ll}T^{1+\frac{A}{4}}\mathcal{L}^{4A},&\text{if }~ 2< A\leq\frac{262}{27},\\T^{\frac{131}{416}A+\frac{154}{416}}\mathcal{L}^{4A+2},&\text{if }~  A>11. \end{array} \right.
\end{align*}
Thus we obtain the estimates for integrals over intervals of the form $[2^{-j}T,2^{1-j}T]$, $1\leq j\ll \log T$, by summing over which, we get Theorem \ref{th1}.
\qed

\section{ Proofs of Theorem \ref{upper bound} and Theorem \ref{one moment}}

Suppose $x\geq1, 1\leq a_i\leq q_i$ and $(a_i, q_i)=1$  for $i=1, 2$. Since periods of both $\cos(2\pi m\theta_1)$ and $\sin(2\pi n\theta_2)$ are 1, we have
 \begin{eqnarray*}
  \nonumber S\big(q_1q_2x; \frac{a_1}{q_1}, \frac{a_2}{q_2}\big)&=&\mathop{{\sum}'}_{mn\leq q_1q_2x} \cos\Big(2\pi m\frac{a_1}{q_1}\Big)\sin\Big(2\pi n\frac{a_2}{q_2}\Big) \\
  \nonumber &=&\sum_{r_1=1}^{q_1}\sum_{r_2=1}^{q_2}\mathop{{\sum}'}_{\begin{subarray}{c} mn\leq q_1q_2x\\m\equiv r_1\!\!\!\!\!\pmod{q_1}\\n\equiv r_2\!\!\!\!\!\pmod{q_2}\end{subarray}} {\cos\Big(2\pi r_1\frac{\displaystyle a_1}{\displaystyle q_1}\Big)\sin\Big(2\pi r_2\frac{\displaystyle a_2}{\displaystyle q_2}\Big)}\\
  &=&\sum_{r_1=1}^{q_1}\sum_{r_2=1}^{q_2}{\cos\Big(2\pi r_1\frac{\displaystyle a_1}{\displaystyle q_1}\Big)\sin\Big(2\pi r_2\frac{\displaystyle a_2}{\displaystyle q_2}\Big)}\mathop{{\sum}'}_{\begin{subarray}{c} mn\leq q_1q_2x\\m\equiv r_1\!\!\!\!\!\pmod{q_1}\\n\equiv r_2\!\!\!\!\!\pmod{q_2}\end{subarray}}{1}.
 \end{eqnarray*}
 Substituting (\ref{4}) into this formula, we obtain
\begin{align}\label{5}
  S\big(q_1q_2x; \frac{a_1}{q_1}, \frac{a_2}{q_2}\big) =&\big(x\log x-x\big)
   \sum_{r_1=1}^{q_{1}}{\cos\Big(2\pi r_1\frac{a_1}{q_1}\Big)}\sum_{r_2=1}^{q_{2}}{\sin\Big(2\pi r_2\frac{a_2}{q_2}\Big)}\\
   \nonumber&-x\sum_{r_1=1}^{q_1}\!{\cos\Big(2\pi r_1\frac{a_1}{q_1}\Big)\frac{\Gamma'}{\Gamma}\Big(\frac{r_1}{q_1}\Big)}
   \sum_{r_2=1}^{q_{2}}{\sin\Big(2\pi r_2\frac{a_2}{q_2}\Big)}\\
 \nonumber&-x\!
 \sum_{r_2=1}^{q_2}{\sin\Big(2\pi r_2\frac{a_2}{q_2}\Big)\frac{\Gamma'}{\Gamma}\Big(\frac{r_2}{q_2}\Big)}\sum_{r_1=1}^{q_{1}}{\cos\Big(2\pi r_1\frac{a_1}{q_1}\Big)}\\
 \nonumber&+\sum_{r_1=1}^{q_1}\sum_{r_2=1}^{q_2}{\cos\big(2\pi r_1\frac{a_1}{ q_1}\big)\sin\big(2\pi r_2\frac{a_2}{q_2}\big)\Delta(q_1q_2x; r_1, q_1, r_2, q_2)}.
\end{align}
 For $(a, q)=1$ and $\theta$ an arbitrary constant, it is easy to see that
 \begin{align}\label{tripro}
    \sum\limits_{r=1}^{q}{\sin(2\pi r\frac{a}{q}+\theta)}=\sum\limits_{r=1}^{q}{\cos(2\pi r\frac{a}{q}+\theta)}=0,
 \end{align}
which implies
\begin{align}\label{8}
  \!S\big(q_1q_2x; \frac{a_1}{q_1}, \frac{a_2}{q_2}\big)\!=\!\!\sum_{r_1=1}^{q_1}\sum_{r_2=1}^{q_2}{\!\cos\big(2\pi\frac{\displaystyle r_1a_1}{\displaystyle q_1}\big)\sin\big(2\pi\frac{\displaystyle r_2a_2}{\displaystyle q_2}\big)\Delta(q_1q_2x; r_1, q_1, r_2, q_2)}.
\end{align}

Using the estimate of $\Delta(x; r_1, q_1, r_2, q_2)$ shown by \eqref{9a}, and M\"{u}ller and Nowak's \eqref{muller}, we complete the proof of Theorem~\ref{upper bound} and Theorem~\ref{one moment}.

\section{ Proofs of Theorem \ref{thmmain2}}

\subsection{ An expression of $S^k(q_1q_2x; \dfrac{a_1}{q_1}, \dfrac{a_2}{q_2})$}

\par
From (\ref{8}), \eqref{lemliu}, noting the fact that $T\gg (q_1q_2)^{\varepsilon'}$, we have
\begin{align*}
  \nonumber &S\big(q_1q_2x; \frac{a_1}{q_1}, \frac{a_2}{q_2}\big)\\
=&R_0(x; y)+R_{12}(x; y, H)+R_{21}(x; y, H)+G_{12}(x; H)+G_{21}(x; H)
  +O\big(q_1q_2\mathcal{L}^3\big).
\end{align*}
Here
\begin{gather}
     R_0(x; y)=\frac {x^{\frac{1}{4}}}{\sqrt{2}\pi}\sum_{r_1=1}^{q_1}\sum_{r_2=1}^{q_2}\cos\big(2\pi r_1\frac{a_1}{q_1}\big)\sin\big(2\pi r_2\frac{a_2}{q_2}\big)\sum_{n\leq y}\frac{\tau(n, x)}{n^{\frac{3}{4}}}, \label{R0} \\
\nonumber R_{12}(x; y, H) = \frac {x^{\frac{1}{4}}}{\sqrt{2}\pi}\sum_{r_1=1}^{q_1}\sum_{r_2=1}^{q_2}\cos\big(2\pi r_1\frac{a_1}{q_1}\big) \sin\big(2\pi
      r_2\frac{a_2}{q_2}\big)\sum_{y<n\leq 2^{J+1}H^2}\frac{\tau_{12}(n, x; H)}{n^{\frac{3}{4}}}, \\
\nonumber  R_{21}(x; y, H) = \frac {x^{\frac{1}{4}}}{\sqrt{2}\pi}\sum_{r_1=1}^{q_1}\sum_{r_2=1}^{q_2}\cos\big(2\pi r_1\frac{a_1}{q_1}\big)\sin\big(2\pi
       r_2\frac{a_2}{q_2}\big)\sum_{y<n\leq 2^{J+1}H^2}\frac{\tau_{21}(n, x; H)}{n^{\frac{3}{4}}},\\
    G_{12}(x; H) = O\Bigl(q_2\sum_{\begin{subarray}{c}  n_1\leq q_1\sqrt{T}\end{subarray}}\min\Bigl(1, \frac{1}{H\|\frac{q_1x}{n_1}-\frac{r_2}{q_2}\|}\Bigr)\Bigr),\label{G12}\\
 \nonumber  G_{21}(x; H)=O\Big(q_1\sum_{\begin{subarray}{c}  n_2\leq q_2\sqrt{T}\end{subarray}}\min\Big(1, \frac{1}{H\|\frac{q_2x}{n_2}-\frac{r_1}{q_1}\|}\Big)\Big).
\end{gather}

Denote ${R}_0=R_0(x; y)$, $R_1=R_{12}(x; y, H)$, $R_2=R_{21}(x; y, H)$, $G_{12}=G_{12}(x; H)$, $G_{21}=G_{21}(x; H)$,
$R=R_0+ R_1+ R_2$, and $G=G_{12}+G_{21}$. Then
\begin{eqnarray}\label{s2}
\nonumber S\big(q_1q_2x; \frac{a_1}{q_1}, \frac{a_2}{q_2}\big)= R+G+O\big(q_1q_2\mathcal{L}^3\big).
\end{eqnarray}
Let $k\geq2$ be a fixed integer. By the elementary formula $(a+b)^k=a^k+O(|a|^{k-1}|b|+|b|^k)$, we get
  \begin{equation}\label{21}
  S^k\big(q_1q_2x; \frac{a_1}{q_1}, \frac{a_2}{q_2}\big)\!=\!R^k+O\big(|R|^{k-1}G+|R|^{k-1}q_1q_2\mathcal{L}^3\big)+O\big(G^k +(q_1q_2\mathcal{L}^{3})^{k}\big).
  \end{equation}
The expression of $R_0$ given by (\ref{R0}) shows
that
\begin{align}\label{R012}
  R_0=&\frac {x^{\frac{1}{4}}}{\sqrt{2}\pi}\sum_{n\leq y}n^{-\frac{3}{4}}\!\sum_{n=hl}\sum_{r_1=1}^{q_1}\sum_{r_2=1}^{q_2}\!\cos\big(2\pi r_1\frac{a_1}{q_1}\big)\sin\big(2\pi r_2\frac{a_2}{q_2}\big)\\
 \nonumber &\times\cos\Big(4\pi \sqrt{nx}-2\pi\Big(\frac{hr_2}{q_2}+\frac{lr_1}{q_1}+\frac{1}{8}\Big)\Big).
\end{align}
  Set
  \begin{equation*}
    S_0=\sum_{r_1=1}^{q_1}\cos\big(2\pi r_1\frac{a_1}{q_1}\big)\cos\Big(4\pi \sqrt{nx}\!-\!2\pi\Big(\frac{hr_2}{q_2}+\frac{lr_1}{q_1}+\frac{1}{8}\Big)\Big).
  \end{equation*}
 Then, clearly
 \begin{multline*}
     S_0=\frac{1}{2}\!\sum_{r_1=1}^{q_1}\!\bigg\{\cos\Big(4\pi \! \sqrt{nx}\!-\!2\pi\Big(\frac{hr_2}{q_2}\!+\!\frac{r_1}{q_1}(l+a_1)\!+\!\frac{1}{8}\Big)\Big)\!\\
     +\!\cos\Big(4\pi \! \sqrt{nx}\!-\!2\pi\Big(\frac{hr_2}{q_2}\!+\!\frac{r_1}{q_1}(l-a_1)\!+\!\frac{1}{8}\Big)\Big)\bigg\}.
 \end{multline*}
 If $l\not\equiv\pm a_1\!\!\!\pmod{q_1}$, from \eqref{tripro}, $S_0=0$.

If $q_1>2$, then for $l\equiv a_1\!\!\!\pmod{q_1}$ or $l\equiv -a_1\!\!\!\pmod{q_1}$, obviously they cannot be simultaneously true,
 \begin{align*}
 S_0=&\frac{1}{2}\sum_{r_1=1}^{q_1}\cos\Big(4\pi \sqrt{nx}-2\pi\big(\frac{hr_2}{q_2}+\frac{1}{8}\big)\Big)=\frac{q_1}{2}\cos\Big(4\pi \sqrt{nx}-2\pi\big(\frac{hr_2}{q_2}+\frac{1}{8}\big)\Big).
 \end{align*}

If $q_1=2$, then $a_1=1$, and $a_1\equiv\!-a_1\!\!\!\pmod{q_1}$. Thus, for $l\equiv a_1(\text{or} -a_1)\!\!\!\pmod{q_1}$,
 \begin{equation*}
 S_0=q_1\cos\Big(4\pi \sqrt{nx}-2\pi\big(\frac{hr_2}{q_2}+\frac{1}{8}\big)\Big).
 \end{equation*}
Summing up, we see
\begin{align}\label{cos1}
   & \sum_{r_1=1}^{q_1}\!\cos\big(2\pi r_1\frac{a_1}{q_1}\big)\cos\Big(4\pi \sqrt{nx}\!-\!2\pi\big(\frac{hr_2}{q_2}+\frac{l r_1}{q_1}+ \frac{1}{8}\big)\Big) \\
     \nonumber=&\left\{\begin{array}
     {ll}
     \dfrac{q_1}{2}\cos\Big(4\pi \sqrt{nx}-2\pi\big(\dfrac{hr_2}{q_2}+\dfrac{1}{8}\big)\Big),&q_1>2, l\equiv a_1,-a_1\!\!\!\!\!\pmod{q_1}; \\
     q_1\cos\Big(4\pi \sqrt{nx}-2\pi\big(\dfrac{hr_2}{q_2}+\dfrac{1}{8}\big)\Big),&q_1=2, l\equiv a_1\!\!\!\!\!\pmod{q_1}; \\
     0,&\text{else}. \end{array}\right.
   \end{align}
By a similar argument, we get
\begin{align}\label{cossin1}
& \sum_{r_2=1}^{q_2}\!\sin\big(2\pi r_2\frac{ a_2}{q_2}\big)\cos\Big(4\pi \sqrt{nx}\!-\!2\pi\big(\frac{hr_2}{q_2}+\frac{1}{8}\big)\Big)\\
  \nonumber =&\left\{\begin{array}
   {ll}
   \dfrac{q_2}{2}\cos\big(4\pi\sqrt{nx}-\dfrac{3\pi}{4}\big),&h\equiv a_2\!\!\!\!\!\pmod{q_2}; \\
   -\dfrac{q_2}{2}\cos\big(4\pi\sqrt{nx}-\dfrac{3\pi}{4}\big),&h\equiv -a_2\!\!\!\!\!\pmod{q_2}; \\
   0,&\text{else}. \end{array}\right.
   \end{align}
Combining (\ref{R012})-(\ref{cossin1}), we obtain
   \begin{align}\label{15}
    R_0=& \frac{x^\frac{1}{4}}{\sqrt{2}\pi}\times\frac{1}{2}\times\frac{1}{2}q_1q_2\sum_{n\leq y}\frac{\cos\big(4\pi\sqrt{nx}-\frac{3\pi}{4}\big)}{n^{\frac{3}{4}}}\\
     \nonumber &\times\bigg(\sum_{\begin{subarray}{c}n=hl\\h\equiv a_2\!\!\!\!\!\pmod{q_2}\\l\equiv a_1\!\!\!\!\!\pmod{q_1}\end{subarray}}1\!+\!\sum_{\begin{subarray}{c}n=hl\\h\equiv a_2\!\!\!\!\!\pmod{q_2}\\l\equiv -a_1\!\!\!\!\!\pmod{q_1}\end{subarray}}1\!-\!\sum_{\begin{subarray}{c}n=hl\\h\equiv -a_2\!\!\!\!\!\pmod{q_2}\\l\equiv a_1\!\!\!\!\!\pmod{q_1}\end{subarray}}1\!-\!\sum_{\begin{subarray}{c}n=hl\\h\equiv -a_2\!\!\!\!\!\pmod{q_2}\\l\equiv -a_1\!\!\!\!\!\pmod{q_1}\end{subarray}}1\bigg) \\
     \nonumber =&\frac{q_1q_2x^\frac{1}{4}}{4\sqrt{2}\pi}\sum_{n\leq y}\frac{\cos\big(4\pi\sqrt{nx}-\frac{3\pi}{4}\big)}{n^{\frac{3}{4}}}\Delta d_2(n;a_1, q_1, a_2, q_2).\hspace{6em}
   \end{align}

Similarly, we have
   \begin{align}\label{16}
     R_1 =\frac{q_1q_2x^\frac{1}{4}}{4\sqrt{2}\pi}\sum_{y<n\leq2^{J+1}H^2 }\frac{\cos\big(4\pi\sqrt{nx}-\dfrac{3\pi}{4}\big)}{n^{\frac{3}{4}}}\Delta d_{2,1}(n,H,J;a_1, q_1, a_2, q_2),\\
   \nonumber R_2 = \frac{q_1q_2x^\frac{1}{4}}{4\sqrt{2}\pi}\sum_{y<n\leq 2^{J+1}H^2 }\frac{\cos\big(4\pi\sqrt{nx}-\dfrac{3\pi}{4}\big)}{n^{\frac{3}{4}}}\Delta d_{2,2}(n,H,J;a_1, q_1, a_2, q_2).
   \end{align}
\subsection{Upper bound of $\int_{\frac{T}{2}}^TG^kdx~~(k\geq1)$}
\par
  From the expression of $G_{12}$ provided by (\ref{G12}), we have
  \begin{align*}
   \nonumber \int_{\frac{T}{2}}^TG_{12}(x; H)dx \ll q_2\sum_{n_1\leq q_1\sqrt{T}}\int_{\frac{T}{2}}^T\min \Big(1, \frac{1}{H\|\frac{q_1x}{n_1}-\frac{r_2}{q_2}\|}\Big)dx.
  \end{align*}
Let $u=\dfrac{q_1x}{n_1}-\dfrac{l_2}{q_2}$. Then
\begin{align}\label{22}
   \int_{\frac{T}{2}}^TG_{12}(x; H)dx\ll& q_2\sum_{n_1\leq q_1\sqrt{T}}\int_{\frac{q_1T}{2n_1}-\frac{r_2}{q_2}}^{\frac{q_1T}{n_1}-\frac{r_2}{q_2}}\min \Big(1, \frac{1}{H\| u\|}\Big)\frac{n_1}{q_1}du \\
   \nonumber\ll& T q_2\sum_{n_1\leq q_1\sqrt{T}}\int_0^1\min \Big(1, \frac{1}{H\| u\|}\Big)du\\
  \nonumber \ll& q_1q_2T^\frac{3}{2}\int_0^\frac{1}{2}\min \Big(1, \frac{1}{Hu}\Big)du.
\end{align}
Since
\begin{align*}
  \nonumber\int_0^\frac{1}{2}\min \Big(1, \frac{1}{Hu}\Big)du=\int_\frac{1}{H}^\frac{1}{2}\frac{1}{Hu}du+\int_0^\frac{1}{H}1du
   \ll \frac{1}{H}\mathcal{L},
\end{align*}
 we obtain
\begin{eqnarray*}
  \int_{\frac{T}{2}}^TG_{12}(x; H)dx\ll q_1 q_2T^\frac{3}{2}H^{-1}\mathcal{L}.
\end{eqnarray*}

 For $\int_{\frac{T}{2}}^TG_{21}(x; H)dx$, this estimate still holds, which implies
\begin{equation*}
  \int_{\frac{T}{2}}^TG(x; H)dx\ll q_1q_2T^\frac{3}{2}H^{-1}\mathcal{L}.
\end{equation*}
Then by $G\ll q_1q_2T^{\frac{1}{2}}$, we get for $k\geq1$,
\begin{align*}
 \nonumber \int_{\frac{T}{2}}^TG^{k}(x)dx \ll& \big(q_1q_2T^{\frac{1}{2}}\big)^{k-1}\int_{\frac{T}{2}}^TG(x)dx\ll (q_1q_2)^{k}T^{1+\frac{k}{2}}H^{-1}\mathcal{L}.
\end{align*}
Taking $H=T^{\frac{A_0}{2}}$, we have the following
\begin{lemma}\label{lem1}{\, \, \rm} If $k\geq1, k\in \mathbb{R}$, then
\begin{equation*}
  \int_{\frac{T}{2}}^TG^{k}(x)dx\ll (q_1q_2)^kT\mathcal{L}.
\end{equation*}
\end{lemma}

\subsection{Evaluation of the integral $\int_{\frac{T}{2}}^TR_0^2(x; y)dx$}

Now we go to evaluate the integral $\int_{\frac{T}{2}}^TR_0^2(x; y)dx$, which provides the main term in the asymptotic formula of $ \int_{\frac{T}{2}}^{T}{S^2(x, \frac{a_1}{q_1}, \frac{a_2}{q_2})}dx$.

 By the expression of $R_0$ given by (\ref{15}), using the elementary formula
 $$\cos a_1\cos a_2=\frac{1}{2}\cos(a_1- a_2)+\frac{1}{2}\cos(a_1+a_2),$$
we get
\begin{equation}\label{s1s2}
  R_0^2(x; y)=\frac{(q_1q_2)^2}{64\pi^2}\big(S_1(x)+S_2(x)\big),
\end{equation}
where
\begin{align*}
    S_1(x)=&x^{\frac{1}{2}}\sum_{n_1,n_2\leq y}\frac{\Delta d_2(n_1;a_1, q_1, a_2, q_2)\Delta d_2(n_2;a_1, q_1, a_2, q_2)}{(n_1n_2)^{\frac{3}{4}}}\\
    &\hspace{9em}\times\cos\big(4\pi\sqrt{x}(\sqrt{n_1}-\sqrt{n_2})\big),\\
    S_2(x)=&x^{\frac{1}{2}}\sum_{n_1,n_2\leq y}\frac{\Delta d_2(n_1;a_1, q_1, a_2, q_2)\Delta d_2(n_2;a_1, q_1, a_2, q_2)}{(n_1n_2)^{\frac{3}{4}}}\\
    &\hspace{9em}\times\cos\big(4\pi\sqrt{x}(\sqrt{n_1}+\sqrt{n_2})-\frac{3\pi}{2}\big).
\end{align*}

  Noticing that for $a, b>1$,
  \begin{eqnarray}\label{33}
    \int_a^b\cos(A\sqrt{t}+B)dt\ll \frac{\sqrt{b}+\sqrt{a}}{|A|},
   \end{eqnarray}
we have
   \begin{align}\label{34}
     \int_{\frac{T}{2}}^Tx^{\frac{k}{4}}\cos(A\sqrt{x}+B)dx=& \int_{\frac{T}{2}}^Tx^{\frac{k}{4}}d\Big(\int_{\frac{T}{2}}^x\cos(A\sqrt{t}+B)dt\Big)\ll T^{\frac{k}{4}+\frac{1}{2}}|A|^{-1}.
   \end{align}
Hence, we get
   \begin{align}\label{ss2}
    \int_{\frac{T}{2}}^{T}S_{2}(x)dx=&\sum_{n_1, n_2\leq y}\frac{\Delta d_2(n_1;a_1, q_1, a_2, q_2)\Delta d_2(n_2;a_1, q_1, a_2, q_2)}{(n_1n_2)^{\frac{3}{4}}}\\
    \nonumber &\times\int_{\frac{T}{2}}^{T} x^{\frac{1}{2}}\cos\big(4\pi\sqrt{x}(\sqrt{n_1} + \sqrt{n_2}) - \frac{3\pi}{2}\big)dx\\
\nonumber\ll&T\sum_{n_1, n_2\leq y}\frac{d(n_1) d(n_2)}{(n_1n_2)^{\frac{3}{4}}(\sqrt{n_1}+\sqrt{n_2})}\\
\nonumber\ll&T\sum_{n_1, n_2\leq y}\frac{d(n_1) d(n_2)}{n_1n_2}\ll T\log^4y,
 \end{align}
by using $\sqrt{n_1}+\sqrt{n_2}\geq2(n_1n_2)^{\frac{1}{4}}$, and $\sum_{n\leq N}\frac{d(n)}{n}\ll\log^2N$ in the last step.

   Now we consider $\int_{\frac{T}{2}}^{T}S_{1}(x)dx$. Denote
   \begin{equation}\label{s10}
    S_1(x)=S_{11}(x)+S_{12}(x),
   \end{equation}
   where
   \begin{align*}
    S_{11}(x)=&x^{\frac{1}{2}}\sum_{n\leq y}\frac{\big(\Delta d_2(n;a_1, q_1, a_2, q_2)\big)^2}{n^{\frac{3}{2}}},\\
    S_{12}(x)=&x^{\frac{1}{2}}\sum_{n_1\neq n_2\leq y}\frac{\Delta d_2(n_1;a_1, q_1, a_2, q_2)\Delta d_2(n_2;a_1, q_1, a_2, q_2)}{(n_1n_2)^{\frac{3}{4}}}\\
    &\hspace{12em}\times\cos\big(4\pi\sqrt{x}(\sqrt{n_1}-\sqrt{n_2})\big).
   \end{align*}
   Then
      \begin{align}\label{s11}
    \int_{\frac{T}{2}}^{T}S_{11}(x)dx=&\sum_{n\leq y}\frac{\big(\Delta d_2(n;a_1, q_1, a_2, q_2)\big)^2}{n^{\frac{3}{2}}}\int_{\frac{T}{2}}^{T}x^{\frac{1}{2}}dx\\
\nonumber=&\sum_{n=1}^{\infty}\frac{\big(\Delta d_2(n;a_1, q_1, a_2, q_2)\big)^2}{n^{\frac{3}{2}}}\int_{\frac{T}{2}}^{T}x^{\frac{1}{2}}dx+O\Big(T^{\frac{3}{2}}\sum_{n> y}\frac{d^2(n)}{n^{\frac{3}{2}}}\Big)\\
\nonumber=&B_2\big(\frac{a_1}{q_1}, \frac{a_2}{q_2}\big)\int_{\frac{T}{2}}^{T}x^{\frac{1}{2}}dx+O(T^{\frac{3}{2}}y^{-\frac{1}{2}}\log^3y),
    \end{align}
    noting that $\sum_{n\leq N}d^2(n)\ll N\log^3 N$.
By using Lemma \ref{lem2.4} and \eqref{34}, we get
   \begin{align}\label{s12}
    \int_{\frac{T}{2}}^{T}S_{12}(x)dx\ll&T\sum_{n_1\neq n_2\leq y}\frac{d(n_1) d(n_2)}{(n_1n_2)^{\frac{3}{4}}|\sqrt{n_1}-\sqrt{n_2}|}\\
\nonumber\ll&T\sum_{n\leq y}\frac{1}{\sqrt{n+1}-\sqrt{n}}\cdot\frac{d^2(n)}{n^\frac{3}{2}}  \\
\nonumber \ll&T\sum_{n\leq y}\frac{d^2(n)}{n}\ll T\log^4y.
    \end{align}
  From \eqref{s10}-\eqref{s12}, we see
 \begin{align}\label{ss1}
 \int_{\frac{T}{2}}^{T}S_{1}(x)dx=&B_2(\frac{a_1}{q_1}, \frac{a_2}{q_2})\int_{\frac{T}{2}}^{T}x^{\frac{1}{2}}dx+O(T^{\frac{3}{2}}y^{-\frac{1}{2}}\log^3y),
    \end{align}
    in view of $y\leq T\mathcal{L}^{-4}$.
Combing \eqref{s1s2}, \eqref{ss2} and \eqref{ss1}, we get

\begin{lemma}\label{lem2.2}
\begin{equation*}
\int_{\frac{T}{2}}^T\!\!R_0^2dx=\frac{(q_1q_2)^2}{2^6\pi^2}B_2\big(\frac{a_1}{q_1}, \frac{a_2}{q_2}\big) \int_{\frac{T}{2}}^Tx^{\frac{1}{2}}dx+O\big((q_1q_2)^2T^{\frac{3}{2}}y^{-\frac{1}{2}}\log^3y\big).
\end{equation*}
\end{lemma}

\subsection{Evaluation of the integral $\int_{\frac{T}{2}}^TR_0(R_1+R_2)dx$}

Denote $N=2^{J+1}H^2$.
From (\ref{15}) and (\ref{16}), we see
\begin{multline}\label{73}
  R_0R_1\!=\!\frac{(q_1q_2)^2x^{\frac{1}{2}}}{2^5\pi^2}\!\!\!\!\sum_{y<n_1\leq N}\sum_{n_2\leq y}\!\frac{\cos\big(4\pi\sqrt{n_1x}-\frac{3\pi}{4}\big)\cos\big(4\pi\sqrt{n_2x}-\frac{3\pi}{4}\big)} {(n_1n_2)^{\frac{3}{4}}}\\
  \times\Delta d_{2,1}(n_1,H,J;a_1, q_1, a_2, q_2)\Delta d_2(n_2;a_1, q_1, a_2, q_2).
\end{multline}
Since
\begin{align*}
  &2\cos\big(4\pi\sqrt{n_1x}-\frac{3\pi}{4}\big)\cos\big(4\pi\sqrt{n_2x}-\frac{3\pi}{4}\big)\\
  =& \cos\big(4\pi(\sqrt{n_1}-\sqrt{n_2})\sqrt{x}\big)-\sin\big(4\pi(\sqrt{n_1}+\sqrt{n_2})\sqrt{x}\big),
\end{align*}
by using (\ref{34}), we get
\begin{align*}
   \nonumber \int_{\frac{T}{2}}^T x^{\frac{1}{2}}\cos\big(4\pi\sqrt{n_1x}-\frac{\pi}{4}\big)\cos\big(4\pi\sqrt{n_2x}-\frac{\pi}{4}\big)dx
     \ll& \frac{T}{\sqrt{n_1}-\sqrt{n_2}}.
  \end{align*}
in view of $n_1\neq n_2$. Thus, we have
\begin{equation*}
  \nonumber \int_{\frac{T}{2}}^T\!\!R_0R_1dx \ll (q_1q_2)^2T\sum_{y<n_1\leq N}\sum_{n_2\leq y}\frac{d(n_1)d(n_2)}{(n_1n_2)^\frac{3}{4}(\sqrt{n_1}-\sqrt{n_2})}.
\end{equation*}
Dividing the interval $y<n_1\leq N$ into $y<n_1\leq 4y$ and $4y<n_1\leq N$, we see
 \begin{eqnarray}\label{74}
   \int_{\frac{T}{2}}^TR_0R_1dx \ll (q_1q_2)^2T\Big(\sum\nolimits_5+\sum\nolimits_6\Big),
 \end{eqnarray}
where
\begin{flalign}
\nonumber\sum\nolimits_5=\sum_{y<n_1\leq 4y}\sum_{n_2\leq y}\frac{d(n_1)d(n_2)}{(n_1n_2)^\frac{3}{4}(\sqrt{n_1}-\sqrt{n_2})}, \\
\nonumber\sum\nolimits_6=\sum_{4y<n_1\leq N}\sum_{n_2\leq y}\frac{d(n_1)d(n_2)}{(n_1n_2)^\frac{3}{4}(\sqrt{n_1}-\sqrt{n_2})}.
\end{flalign}
By the same agrement with \eqref{s12}, we get
\begin{eqnarray}\label{sum5}
\sum\nolimits_5&\ll&\!\!\!\!\sum_{\begin{subarray}{c}n_1, n_2\leq 4y\\n_1\neq n_2\end{subarray}}\frac{d(n_1)d(n_2)}{(n_1n_2)^\frac{3}{4}(\sqrt{n_1}-\sqrt{n_2})}\ll T\log^4y.
\end{eqnarray}
It is clear that
\begin{align}\label{sum6}
\sum\nolimits_6\ll \sum_{4y<n_1\leq N}\frac{d(n_1)}{n_1^\frac{5}{4}}\sum_{n_2\leq y}\frac{d(n_2)}{n_2^\frac{3}{4}}\ll \log^2y.
\end{align}
Combining (\ref{74})-\eqref{sum6}, we have
\begin{equation*}
  \int_{\frac{T}{2}}^TR_0R_1dx \ll (q_1q_2)^2T\log^4y.
\end{equation*}
The same estimate holds for $\int_{\frac{T}{2}}^TR_0R_2dx$, which yields

\begin{lemma}\label{lem3.2}
\begin{equation*}
\int_{\frac{T}{2}}^TR_0(R_1+R_2)dx \ll (q_1q_2)^2T\log^4y.
\end{equation*}
\end{lemma}

\subsection{ Mean-square of $R_1+R_2$}

From the expression of $R_1$ given by (\ref{16}), we get
\begin{equation*}
  R_1\ll x^{\frac{1}{4}}q_1q_2\Big|\sum_{y<n\leq N}\frac{e(2\sqrt{nx}-\frac{3}{8})}{n^{\frac{3}{4}}}\Delta d_{2,1}(n,H,J;a_1, q_1, a_2, q_2)\Big|.
\end{equation*}
So we have
\begin{align*}
  R_1^2 \ll&x^{\frac{1}{2}}(q_1q_2)^2\Big|\sum_{y<n\leq N}\frac{e(2\sqrt{nx}-\frac{3}{8})}{n^{\frac{3}{4}}}\Delta d_{2,1}(n,H,J;a_1, q_1, a_2, q_2)\Big|^2  \\
=& x^{\frac{1}{2}}(q_1q_2)^2\!\!\!\sum_{y<n_1,n_2\leq N}\!\!\frac{e\big(2(\sqrt{n_2}-\sqrt{n_1})\sqrt{x}\big)}{(n_1n_2)^{\frac{3}{4}}}\prod_{i=1}^2 \Delta d_{2,1}(n_i,H,J;a_1, q_1, a_2, q_2).
\end{align*}
Splitting the sum  according to $n_1=n_2$ or $n_1\neq n_2$, we obtain
\begin{multline*}
   R_1^2\ll x^{\frac{1}{2}}(q_1q_2)^2\Biggl\{\sum_{y<n\leq N}\frac{d^2(n)}{n^{\frac{3}{2}}}+\!\!\sum_{y<n_1<n_2\leq N}\frac{e\big(2(\sqrt{n_2}-\sqrt{n_1})\sqrt{x}\big)}{(n_1n_2)^{\frac{3}{4}}}\\
   \times\prod_{i=1}^2 d_{2,1}(n_i,H,J;a_1, q_1, a_2, q_2) \Biggr\}.
\end{multline*}
Similarly to (\ref{34}), we get the estimate
\begin{equation*}
  \int_{\frac{T}{2}}^Tx^\frac{1}{2}e\big(2(\sqrt{n_2}-\sqrt{n_1})\sqrt{x}\big)dx\ll\frac{T}{\sqrt{n_2}-\sqrt{n_1}},
\end{equation*}
which yields
\begin{align}\label{14}
\int_{\frac{T}{2}}^TR_1^2dx\ll& (q_1q_2)^2\! \sum_{y<n\leq N}\frac{d^2(n)}{n^{\frac{3}{2}}}\int_{\frac{T}{2}}^Tx^{\frac{1}{2}}dx+\!\!(q_1q_2)^2\\
\nonumber&\times\sum_{y<n_1<n_2\leq N}\frac{d(n_1)d(n_2)}{(n_1n_2)^{\frac{3}{4}}}\int_{\frac{T}{2}}^Tx^{\frac{1}{2}}e\big(2(\sqrt{n_2}-\sqrt{n_1})\sqrt{x}\big)
   dx\\
\nonumber\ll& T^{\frac{3}{2}}(q_1q_2)^2y^{-\frac{1}{2}}\log^3y\!+\!T(q_1q_2)^2\!\!\!\!\sum_{y<n_1<n_2\leq N}\!\frac{d(n_1)d(n_2)}{(n_1n_2)^{\frac{3}{4}}(\sqrt{n_2}-\sqrt{n_1})},
\end{align}
where we used the estimate $\sum\limits_{n\leq u}d^2(n)\ll u\log^3u$.
By Lemma \ref{lem2.4}, we have
\begin{align*}
 \nonumber \sum_{y<n_1<n_2\leq N}\frac{d(n_1)d(n_2)}{(n_1n_2)^{\frac{3}{4}}(\sqrt{n_2}-\sqrt{n_1})}\ll \sum_{y<n\leq N} \frac{d^2(n)}{n}
  \ll \mathcal{L}^{4} .
\end{align*}
From \eqref{14}, we obtain the estimate
\begin{align*}
 \nonumber\int_{\frac{T}{2}}^TR_1^2dx
 \ll& (q_1q_2)^2T^{\frac{3}{2}}y^{-\frac{1}{2}}\mathcal{L}^{3}+(q_1q_2)^2 T\mathcal{L}^{4}.
\end{align*}

For $\int_{\frac{T}{2}}^TR_2^2dx$, the same estimate holds. Thus, we get
\begin{lemma}\label{lem4.2}
\begin{equation*}
\int_{\frac{T}{2}}^T|R_1+R_2|^2dx\ll (q_1q_2)^2T^{\frac{3}{2}}y^{-\frac{1}{2}}\mathcal{L}^{3}+(q_1q_2)^2 T\mathcal{L}^{4}.
\end{equation*}
\end{lemma}

\subsection{Evaluation of the integral $\int_{\frac{T}{2}}^TR^2dx$}　

Since $R=R_0+R_1+R_2$, we have
\begin{equation*}
  \int_{\frac{T}{2}}^T\!\!R^2dx=\!\!\int_{\frac{T}{2}}^T\!\!R_0^2dx+2\int_{\frac{T}{2}}^T\!\!R_0(R_1+R_2)dx+ \int_{\frac{T}{2}}^T(R_1+R_2)^2dx.
\end{equation*}
Then from Lemma \ref{lem2.2}, Lemma \ref{lem3.2} and Lemma \ref{lem4.2}, it follows that
 \begin{equation*}
\int_{\frac{T}{2}}^TR^2dx = \frac{(q_1q_2)^2}{2^6\pi^2}B_2\big(\frac{a_1}{q_1}, \frac{a_2}{q_2}\big) \int_{\frac{T}{2}}^Tx^{\frac{1}{2}}dx+O\big((q_1q_2)^2T^{\frac{3}{2}}y^{-\frac{1}{2}}\mathcal{L}^{3}\big).
\end{equation*}
Take $y=T^\frac{1}{2}$. We see
\begin{lemma}\label{lem5.2}
\begin{equation*}
\int_{\frac{T}{2}}^TR^2dx=
\frac{(q_1q_2)^2}{2^6\pi^2}B_2\big(\frac{a_1}{q_1}, \frac{a_2}{q_2}\big) \int_{\frac{T}{2}}^Tx^{\frac{1}{2}}dx+O\big((q_1q_2)^2T^{\frac{5}{4}}\mathcal{L}^{3}\big).
\end{equation*}

\end{lemma}

\subsection{Mean-value of $S(q_1q_2x, \dfrac{a_1}{q_1}, \dfrac{a_2}{q_2})$}　

From Lemma~\ref{lem5.2} and Lemma~\ref{lem1}, using Cauchy-Schwarz inequality, we get
\begin{align}\label{61.2}
\nonumber\int_{\frac{T}{2}}^T|R|Gdx &\ll(q_1q_2)^2T^{\frac{5}{4}}\mathcal{L}^{2},\\
  \nonumber q_1q_2\mathcal{L}^{3}\int_{\frac{T}{2}}^T|R|dx&\ll(q_1q_2)^{2}T^{\frac{5}{4}}\mathcal{L}^{3},\\
   \nonumber\int_{\frac{T}{2}}^T(q_1q_2)^2\mathcal{L}^6dx&\ll (q_1q_2)^{2}T\mathcal{L}^6.
\end{align}
 Combining these three estimates with Lemma \ref{lem5.2}, Lemma~\ref{lem1}, we obtain
\begin{equation*}
\int_{\frac{T}{2}}^TS^2(q_1q_2x, \dfrac{a_1}{q_1}, \dfrac{a_2}{q_2})dx=
\frac{(q_1q_2)^2}{2^6\pi^2}B_2\big(\frac{a_1}{q_1}, \frac{a_2}{q_2}\big) \int_{\frac{T}{2}}^Tx^{\frac{1}{2}}dx+O\big((q_1q_2)^2T^{\frac{5}{4}}\mathcal{L}^3\big),
\end{equation*}
which immediatly implies Theorem \ref{thmmain2} by noting that the interval $[1,T]$ can be divided into subintervals of the the form $[2^{-j}T,2^{1-j}T]$, $1\leq j\ll \log T$.

\section{ Proofs of Theorem \ref{thmmaink}}

In this section, we shall refer to a method originated from Zhai \cite{zhai04} to prove Theorem \ref{thmmaink}. Throughout this section, we assume $k$ is an integer with $3\leq k<A_0$, where $A_{0}>9$ satisfies \eqref{pre}.
\subsection{Evaluation of the integral $\int_{\frac{T}{2}}^TR_0^k(x; y)dx$}

In this subsection, We will evaluate the integral $\int_{\frac{T}{2}}^TR_0^k(x; y)dx$, which provides the main term in the asymptotic formula of $  \int_{\frac{T}{2}}^{T}{S^k(x, \frac{a_1}{q_1}, \frac{a_2}{q_2})}dx$.

By the elementary formula
\begin{equation}\label{26}
  (\cos a_1)(\cos a_2)\cdots(\cos a_k)\!=\!\frac{1}{2^{k-1}}\!\!\sum_{\textbf{\textit{i}}\in \textbf{1}^{k-1}}\cos \big(a_1+(-1)^{i_1}a_2+\cdots+(-1)^{i_{k-1}}a_k\big),
\end{equation}
from the expression of $R_0$ shown by (\ref{15}), we get
\begin{align*}
R_0^k(x; y)
=&\Big(\frac{x^{\frac{1}{4}}q_1q_2}{4\sqrt{2}\pi}\Big)^k\!\!\frac{1}{2^{k-1}}\sum_{\textbf{\textit{i}}\in \textbf{1}^{k-1}}\\
&\times\sum_{n_i\leq y,1\leq i\leq k}\!\!\!\!\!\frac{\cos\big(4\pi \alpha(\textbf{\textit{n}}; \textbf{\textit{i}})\sqrt{x}\!-\!\beta(\textbf{\textit{i}})\big)\prod_{j=1}^k\!\Delta d_2(n_j;a_1, q_1, a_2, q_2)}{(n_1n_2\cdots n_k)^{\frac{3}{4}}}.
\end{align*}
Denote
\begin{align}
  S_5(x) &= x^{\frac{k}{4}}\sum_{\textbf{\textit{i}}\in \textbf{1}^{k-1}}\sum_{\begin{subarray}{c}n_i\leq y\\1\leq i\leq k\\\alpha(\textbf{\textit{n}}; \textbf{\textit{i}})=0\end{subarray}}\frac{\cos\big(\beta(\textbf{\textit{i}})
  \big)\prod_{j=1}^k\Delta d_2(n_j;a_1, q_1, a_2, q_2)}{(n_1n_2\cdots n_k)^{\frac{3}{4}}}, \label{29}\\
  S_6(x) &= x^{\frac{k}{4}}\!\!\!\sum_{\textbf{\textit{i}}\in \textbf{1}^{k-1}}\!\!\!\sum_{\begin{subarray}{c}n_i\leq y\\1\leq i\leq k\\\alpha(\textbf{\textit{n}}; \textbf{\textit{i}})\neq0\end{subarray}}\!\!\!\!\frac{\cos\big(4\pi \alpha(\textbf{\textit{n}}; \textbf{\textit{i}})\sqrt{x}\!-\!\beta(\textbf{\textit{i}})\big)\prod_{j=1}^k\Delta d_2(n_j;a_1, q_1, a_2, q_2)}{(n_1n_2\cdots n_k)^{\frac{3}{4}}}.\label{30}
\end{align}
Then $R_0^k(x; y)$ can be rewritten as
\begin{equation}\label{28}
  R_0^k(x; y)=\frac{(q_1q_2)^k}{2^{\frac{7}{2}k-1}\pi^k}\big(S_5(x)+S_6(x)\big).
\end{equation}

First we consider $\int_{\frac{T}{2}}^TS_5(x)dx$. By Lemma \ref{lem2.3}, we have
  \begin{align*}
    &\int_{\frac{T}{2}}^TS_5(x)dx
    =\sum_{\textbf{\textit{i}}\in \textbf{1}^{k-1}}\sum_{\begin{subarray}{c}n_i\leq y\\1\leq i\leq k\\\alpha(\textbf{\textit{n}}; \textbf{\textit{i}})=0\end{subarray}}\frac{\cos\beta(\textbf{\textit{i}})}
  {(n_1n_2\cdots n_k)^{\frac{3}{4}}}\prod_{j=1}^k\Delta d_2(n_j;a_1, q_1, a_2, q_2)\int_{\frac{T}{2}}^Tx^{\frac{k}{4}}dx \\
     =& \sum_{\textbf{\textit{i}}\in \textbf{1}^{k-1}}\sum_{\alpha(\textbf{\textit{n}}; \textbf{\textit{i}})=0}\cos\beta(\textbf{\textit{i}})
  \prod_{j=1}^k\frac{\Delta d_2(n_j;a_1, q_1, a_2, q_2)}{n_j^{\frac{3}{4}}} \int_{\frac{T}{2}}^Tx^{\frac{k}{4}}dx\!+\!O(T^{1+\frac{k}{4}}y^{-\frac{1}{2}+\varepsilon}).
  \end{align*}
 Note that
  \begin{align}\label{322}
    \nonumber &\sum\limits_{\textbf{\textit{i}}\in \textbf{1}^{k-1}}\sum\limits_{\alpha(\textbf{\textit{n}}; \textbf{\textit{i}})=0}\cos\beta(\textbf{\textit{i}})
  \prod\limits_{j=1}^k\dfrac{\Delta d_2(n_j;a_1, q_1, a_2, q_2)}{n_j^{\frac{3}{4}}}\\
  \nonumber=&\sum_{\textbf{\textit{i}}\in \textbf{1}^{k-1}}\cos\big((k-2|\textbf{\textit{i}}|)\frac{3\pi}{4}\big)\sum_{\alpha(\textbf{\textit{n}}; \textbf{\textit{i}})=0}\prod_{j=1}^k\frac{\Delta d_2(n_j;a_1, q_1, a_2, q_2)}
  {n_j^{\frac{3}{4}}}  \\
   \nonumber\!=& \sum_{\upsilon=1}^{k-1}\!\cos\big((k\!-\!2\upsilon)\frac{3\pi}{4}\big){k-1\choose \upsilon}\!\!\!\sum_{\sqrt{n_1}+\cdots+\sqrt{n_\upsilon}=\sqrt{n_{\upsilon+1}}+\cdots+\sqrt{n_k}}\!\prod_{j=1}^{k} \frac{\Delta d_2(n_j;a_1, q_1, a_2, q_2)}{n_j^\frac{3}{4}}.
  \end{align}
  According to the definition in Section \ref{yubei}, this is  $B_k\big(\frac{a_1}{q_1}, \frac{a_2}{q_2}\big)$. Thus
\begin{equation}\label{321}
\int_{\frac{T}{2}}^TS_5(x)dx=B_k\big(\frac{a_1}{q_1}, \frac{a_2}{q_2}\big)\int_{\frac{T}{2}}^Tx^{\frac{k}{4}}dx +O(T^{1+\frac{k}{4}}y^{-\frac{1}{2}+\varepsilon}).
\end{equation}

  We now consider $\int_{\frac{T}{2}}^TS_6(x)dx$.
From (\ref{30}), using \eqref{34}, we see
   \begin{align}\label{ints2}
     &\int_{\frac{T}{2}}^TS_6(x)dx\\
     \nonumber =&\!\sum_{\textbf{\textit{i}}\in \textbf{1}^{k-1}}\sum_{\begin{subarray}{c}n_i\leq y,1\leq i\leq k\\\alpha(\textbf{\textit{n}}; \textbf{\textit{i}})\neq0\end{subarray}}\prod_{j=1}^{k}\frac{\Delta d_2(n_j;a_1, q_1, a_2, q_2)}{n_j^\frac{3}{4}} \!\int_{\frac{T}{2}}^T\!\!x^{\frac{k}{4}}\!\cos\big(4\pi \alpha(\textbf{\textit{n}}; \textbf{\textit{i}})\sqrt{x}\!-\!\beta(\textbf{\textit{i}})\big)dx\\
    \nonumber  \ll&T^{\frac{1}{2}+\frac{k}{4}}\sum_{\textbf{\textit{i}}\in \textup{\textbf{1}}^{k-1}}\sum_{\begin{subarray}{c}n_i\leq y,1\leq i\leq k\\\alpha(\textbf{\textit{n}}; \textbf{\textit{i}})\neq0\end{subarray}}\frac{d(n_1)\cdots d(n_k)}{(n_1\cdots n_k)^\frac{3}{4} |\alpha(\textbf{\textit{n}};\textbf{\textit{i}})|}.
   \end{align}
   By Lemma~\ref{lemsum}, we have the estimate
   \begin{eqnarray*}
     \nonumber \sum_{\textbf{\textit{i}}\in \textup{\textbf{1}}^{k-1}}\sum_{\begin{subarray}{c}n_i\leq y,1\leq i\leq k\\\alpha(\textbf{\textit{n}}; \textbf{\textit{i}})\neq0\end{subarray}}\frac{d(n_1)\cdots d(n_k)}{(n_1\cdots n_k)^\frac{3}{4} \alpha(\textbf{\textit{n}};\textbf{\textit{i}})}
\ll y^{s(k)+\varepsilon}.
   \end{eqnarray*}
   Substituting it into (\ref{ints2}), we get
   \begin{equation}
    \int_{\frac{T}{2}}^TS_6(x)dx\ll T^{\frac{1}{2}+\frac{k}{4}}y^{s(k)+\varepsilon}.
   \end{equation}

 Combining this with (\ref{28}), (\ref{321}),  we arrive at
\begin{lemma}\label{lem2}
\begin{equation*}
\int_{\frac{T}{2}}^T\!\!R_0^kdx= \frac{(q_1q_2)^k}{2^{\frac{7}{2}k-1}\pi^k}B_k\big(\frac{a_1}{q_1}, \frac{a_2}{q_2}\big) \int_{\frac{T}{2}}^Tx^{\frac{k}{4}}dx+O\Big((q_1q_2)^k T^{\frac{1}{2}+\frac{k}{4}}\big(T^\frac{1}{2}y^{-\frac{1}{2}+\varepsilon}+y^{s(k)+\varepsilon}\big) \Big).
\end{equation*}
\end{lemma}

\subsection{Evaluation of the integral $\int_{\frac{T}{2}}^TR_0^{k-1}(R_1+R_2)dx$}

Let $N=2^{J+1}H^2$. From the expressions of $R_0$ and $R_1$ given by (\ref{15}) and (\ref{16}), using (\ref{26}), we have
\begin{align*}
 R_0^{k-1}R_1=\big(\frac{q_1q_2x^\frac{1}{4}}{4\sqrt{2}\pi}\big)^k\sum_{y<n_1\leq N}\!\!\sum_{\begin{subarray}{c}n_i\leq y\\2\leq i\leq k\end{subarray}}\!\frac{\cos(4\pi\sqrt{n_1x}-\frac{3\pi}{4})\cdots\cos(4\pi\sqrt{n_kx}-\frac{3\pi}{4})} {(n_1n_2\cdots n_k)^{\frac{3}{4}}}&\\
 \nonumber\times\Delta d_{2,1}(n_1,H,J;a_1, q_1, a_2, q_2)\prod_{i=2}^k\Delta d_2(n_i;a_1, q_1, a_2, q_2)&\\
  \nonumber = (\frac{q_1q_2x^\frac{1}{4}}{4\sqrt{2}\pi})^k\frac{1}{2^{k-1}}\!\!\sum_{\textbf{\textit{i}}\in \textbf{1}^{k-1}}\sum_{y<n_1\leq N}\sum_{\begin{subarray}{c}n_i\leq y\\2\leq i\leq k\end{subarray}}\frac{\cos\big(4\pi \alpha(\textbf{\textit{n}}; \textbf{\textit{i}})\sqrt{x}-\beta(\textbf{\textit{i}})\big)} {(n_1n_2\cdots n_k)^{\frac{3}{4}}}\ \ \ ~&\\
  \nonumber\times\Delta d_{2,1}(n_1,H,J;a_1, q_1, a_2, q_2)\prod_{i=2}^k\Delta d_2(n_i;a_1, q_1, a_2, q_2).&
\end{align*}

Set
\begin{align*}
S_3(x)=&x^{\frac{k}{4}}\sum_{\textbf{\textit{i}}\in \textbf{1}^{k-1}}\sum_{y<n_1\leq N}\sum_{\begin{subarray}{c}n_i\leq y\\2\leq i\leq k\\\alpha(\textbf{\textit{n}}; \textbf{\textit{i}})=0\end{subarray}}\frac{\cos\big(\beta(\textbf{\textit{i}})
  \big)}{(n_1n_2\cdots n_k)^{\frac{3}{4}}}\\
  &\times\Delta d_{2,1}(n_1,H,J;a_1, q_1, a_2, q_2)\prod_{i=2}^k\Delta d_2(n_i;a_1, q_1, a_2, q_2), \\
  S_4(x)=&x^{\frac{k}{4}}\sum_{\textbf{\textit{i}}\in \textbf{1}^{k-1}}\sum_{y<n_1\leq N}\sum_{\begin{subarray}{c}n_i\leq y\\2\leq i\leq k\\\alpha(\textbf{\textit{n}}; \textbf{\textit{i}})\neq0\end{subarray}}\frac{\cos\big(4\pi \alpha(\textbf{\textit{n}}; \textbf{\textit{i}})\sqrt{x}-\beta(\textbf{\textit{i}})\big)} {(n_1n_2\cdots n_k)^{\frac{3}{4}}}\\
  &\times\Delta d_{2,1}(n_1,H,J;a_1, q_1, a_2, q_2)\prod_{i=2}^k\Delta d_2(n_i;a_1, q_1, a_2, q_2).
\end{align*}
Then
\begin{equation}\label{37}
  R_0^{k-1}R_1=\frac{(q_1q_2)^k}{2^{\frac{7}{2}k-1}\pi^k}\big(S_3(x)+S_4(x)\big).
\end{equation}

From Lemma~\ref{lem2.3}, it is easy to get
\begin{align}\label{40}
 \int_{\frac{T}{2}}^TS_3(x)dx&\ll\sum_{\textbf{\textit{i}}\in \textbf{1}^{k-1}}\sum_{y<n_1\leq N}\sum_{\begin{subarray}{c}n_i\leq y\\2\leq i\leq k\\\alpha(\textbf{\textit{n}}; \textbf{\textit{i}})=0\end{subarray}}\prod_{j=1}^{k}\frac{d(n_j)}{n_j^\frac{3}{4}}\int_{\frac{T}{2}}^Tx^{\frac{k}{4}}dx  \\
  \nonumber &\ll T^{1+\frac{k}{4}}\sum_{\upsilon=1}^{k-1}\big|s_{k, \upsilon}(d)-s_{k, \upsilon}(d; y)\big| \ll T^{1+\frac{k}{4}}y^{-\frac{1}{2}+\varepsilon}.
 \end{align}

 We now evaluate $\int_{\frac{T}{2}}^TS_4(x)dx$. From (\ref{34}), we have
 \begin{equation}\label{41}
   \int_{\frac{T}{2}}^TS_4(x)dx\ll  T^{\frac{1}{2}+\frac{k}{4}}\Big(\sum\nolimits_1+\sum\nolimits_2\Big),
 \end{equation}
 where
 \begin{flalign}
 \nonumber&\sum\nolimits_1=\sum_{\textbf{\textit{i}}\in \textbf{1}^{k-1}}\sum_{y<n_1\leq k^2y}\sum_{\begin{subarray}{c}n_i\leq y\\2\leq i\leq k\\\alpha(\textbf{\textit{n}}; \textbf{\textit{i}})\neq0\end{subarray}}\frac{d(n_1)\cdots d(n_k)} {(n_1n_2\cdots n_k)^{\frac{3}{4}}|\alpha(\textbf{\textit{n}}; \textbf{\textit{i}})|}, &&\\
 \nonumber&\sum\nolimits_2=\sum_{\textbf{\textit{i}}\in \textbf{1}^{k-1}}\sum_{ k^2y<n_1\leq N}\sum_{\begin{subarray}{c}n_i\leq y\\2\leq i\leq k\\\alpha(\textbf{\textit{n}}; \textbf{\textit{i}})\neq0\end{subarray}}\frac{d(n_1)\cdots d(n_k)} {(n_1n_2\cdots n_k)^{\frac{3}{4}}|\alpha(\textbf{\textit{n}}; \textbf{\textit{i}})|}. &&
\end{flalign}
 By Lemma \ref{lemsum}, we see
 \begin{equation}\label{42}
   \sum\nolimits_1\ll y^{s(k)+\varepsilon}.
 \end{equation}
For $n_1>k^2y$, obviously $|\alpha(\textbf{\textit{n}} ; \textbf{\textit{i}})|\gg n_1^{\frac{1}{2}}$ holds. Thus,
\begin{align}\label{43}
  \sum\nolimits_2 \ll \sum_{ k^2y<n_1\leq N}\sum_{\begin{subarray}{c}n_i\leq y\\2\leq i\leq k\end{subarray}}\frac{d(n_1)\cdots d(n_k)} {(n_2\cdots n_k)^{\frac{3}{4}}n_1^{\frac{5}{4}}} \ll y^{\frac{k-2}{4}}\mathcal{L}^k.
\end{align}
For $k\geq3$, we see $s(k)>\frac{k-2}{4}$.
Hence, by (\ref{41})-(\ref{43}), we get
\begin{eqnarray}\label{44}
  \int_{\frac{T}{2}}^TS_4(x)dx
   \ll T^{\frac{1}{2}+\frac{k}{4}}y^{s(k)+\varepsilon}.
\end{eqnarray}

From (\ref{37}), (\ref{40}) and (\ref{44}), it follows that
\begin{equation*}
  \int_{\frac{T}{2}}^TR_0^{k-1}R_1dx\ll\frac{(q_1q_2)^k}{2^{\frac{7}{2}k-1}\pi^k} \big(T^{1+\frac{k}{4}}y^{-\frac{1}{2}+\varepsilon}+T^{\frac{1}{2}+\frac{k}{4}}y^{s(k)+\varepsilon}\big).
\end{equation*}

Replacing $R_1$ by $R_2$, this estimate still holds, which immediately implies
\begin{lemma}\label{lem3}
\begin{equation*}
  \int_{\frac{T}{2}}^TR_0^{k-1}(R_1+R_2)dx\ll(q_1q_2)^k \big(T^{1+\frac{k}{4}}y^{-\frac{1}{2}+\varepsilon}+T^{\frac{1}{2}+\frac{k}{4}}y^{s(k)+\varepsilon}\big).
\end{equation*}
\end{lemma}

\subsection{Moments of $R_1+R_2$}

Note $K_0\geq A_0$ is an even number. Let $y^{2s(K_0)}\leq T$. Then from Lemma~\ref{lem2}, we see
\begin{eqnarray*}
  \int_{\frac{T}{2}}^T|R_0|^{K_0}dx \ll (q_1q_2)^{K_0}T^{1+\frac{K_0}{4}+\varepsilon}.
\end{eqnarray*}
By Lemma \ref{lem2.2}, we get
\begin{equation*}
  \int_{\frac{T}{2}}^T\!|R_0|^2dx \ll (q_1q_2)^2T^{\frac{3}{2}+\varepsilon}.
\end{equation*}
Using H\"{o}lder's inequality, we have
\begin{align}\label{51}
  \int_{\frac{T}{2}}^T\!\!|R_0|^{A_0}dx \ll (q_1q_2)^{A_0}T^{1+\frac{A_0}{4}+\varepsilon}.
\end{align}
Owing to (\ref{pre}) and (\ref{8}), we see
\begin{equation}\label{54}
  \int_{\frac{T}{2}}^T\big|S\big(q_1q_2x; \frac{a_1}{q_1}, \frac{a_2}{q_2}\big)\big|^{A_0}dx\ll (q_1q_2)^{A_0}T^{1+\frac{A_0}{4}+\varepsilon}.
\end{equation}
Thus, from Lemma~\ref{lem1}, (\ref{s2}), (\ref{51}) and(\ref{54}), we have
\begin{align}\label{55}
\int_{\frac{T}{2}}^T|R_1+R_2|^{A_0}dx \ll&\int_{\frac{T}{2}}^T\big|S\big(q_1q_2x; \frac{a_1}{q_1}, \frac{a_2}{q_2}\big)\big|^{A_0}dx+\! \int_{\frac{T}{2}}^T|R_0|^{A_0}dx\\
 \nonumber&+ \int_{\frac{T}{2}}^TG^{A_0}(x)dx+O\big(T(q_1q_2)^{A_0}\mathcal{L}^{3A_0}\big)\\
   \nonumber\ll&(q_1q_2)^{A_0}T^{1+\frac{A_0}{4}+\varepsilon}.
\end{align}
 By using Lemma 6.4 and (7.16), and H\"{o}lder's inequality,, we obtain
\begin{lemma}\label{lem4} For $T^\varepsilon\leq y\leq T^{\frac{1}{2s(K_0)}}$, we have
\begin{equation*}
\int_{\frac{T}{2}}^T|R_1+R_2|^{k}dx\ll
(q_1q_2)^{k}T^{1+\frac{k}{4}+\varepsilon}y^{-\frac{A_0-k}{2(A_0-2)}}.
\end{equation*}
\end{lemma}

\subsection{Evaluation of the integral $\int_{\frac{T}{2}}^TR^kdx$}　

Since $(a+b)^k=a^k+ka^{k-1}b+O(|a|^{k-2}|b|^2+|b|^k)$,  we have
\begin{multline*}
  \int_{\frac{T}{2}}^TR^kdx=\int_{\frac{T}{2}}^TR_0^kdx+k\int_{\frac{T}{2}}^TR_0^{k-1}(R_1+R_2)dx\\ +O\Big(\int_{\frac{T}{2}}^T\big|R_0^{k-2}(R_1+R_2)^2\big|dx+\int_{\frac{T}{2}}^T\big|R_1+R_2\big|^kdx\Big).
\end{multline*}
From (\ref{51}), Lemma \ref{lem4} and H\"{o}lder's inequality, we get
\begin{align*}
  \nonumber\int_{\frac{T}{2}}^T\big|R_0^{k-2}(R_1+R_2)^2\big|dx
   \ll (q_1q_2)^{k}T^{1+\frac{k}{4}+\varepsilon}y^{-\frac{A_0-k}{2(A_0-2)}}.
\end{align*}

Combining this two formulas with Lemma \ref{lem2}, Lemma \ref{lem3} and Lemma \ref{lem4}, we obtain
\begin{align*}
   \int_{\frac{T}{2}}^TR^kdx =&\frac{(q_1q_2)^k}{2^{\frac{7}{2}k-1}\pi^k} B_k(\frac{a_1}{q_1}, \frac{a_2}{q_2})\int_{\frac{T}{2}}^Tx^{\frac{k}{4}}dx +O\big((q_1q_2)^{3k}T^{1+\frac{k}{4}+\varepsilon}y^{-\frac{A_0-k}{2(A_0-2)}}\big)\\
   &+O\big((q_1q_2)^kT^{1+\frac{k}{4}}y^{-\frac{1}{2}+\varepsilon}+(q_1q_2)^kT^{\frac{1}{2}+\frac{k}{4}} y^{s(k)+\varepsilon}\big).
\end{align*}
Take $y=T^\frac{1}{2s(K_0)}$,  then $y^{s(k)+\frac{1}{2}}\ll y^{s(K_0)}=T^{\frac{1}{2}}$ and thereby $y^{s(k)}\ll T^{\frac{1}{2}}y^{-\frac{1}{2}}$. Thus, we have the asymptotic formula of $\int_{\frac{T}{2}}^TR^kdx$ as following
\begin{lemma}\label{lem5}
\begin{equation*}
\int_{\frac{T}{2}}^TR^kdx=\frac{(q_1q_2)^k}{2^{\frac{7}{2}k-1}\pi^k}B_k(\frac{a_1}{q_1}, \frac{a_2}{q_2}) \int_{\frac{T}{2}}^Tx^{\frac{k}{4}}dx+O\big((q_1q_2)^{k} T^{1+\frac{k}{4}-\frac{A_0-k}{4(A_0-2)s(K_0)}+\varepsilon}\big).
\end{equation*}
\end{lemma}

\subsection{Moments of $S(q_1q_2x; \dfrac{a_1}{q_1}, \dfrac{a_2}{q_2})$}　

From (\ref{51}) and Lemma~\ref{lem1}, by H\"{o}lder's inequality, we get
\begin{align}
  \label{61}\int_{\frac{T}{2}}^T|R|^{k-1}Gdx \ll(q_1q_2)^kT^{\frac{3}{4}+\frac{k}{4}+\varepsilon},\\
 \label{62}q_1q_2\mathcal{L}^{3}\int_{\frac{T}{2}}^T|R|^{k-1}dx\ll (q_1q_2)^{k}T^{\frac{3}{4}+\frac{k}{4}+\varepsilon}.
\end{align}
Note that
\begin{equation}\label{63}
  \int_{\frac{T}{2}}^T(q_1q_2)^{k}\mathcal{L}^{3k}dx\ll (q_1q_2)^{k}T^{1+\varepsilon}.
\end{equation}
 Since $\frac{A_0-k}{4(A_0-2)s(K_0)}<\dfrac{1}{4}$ for $3\leq k<A_0$, from (\ref{21}), Lemma \ref{lem5}, Lemma~\ref{lem1}, (\ref{61})-(\ref{63}), we can conclude that
\begin{multline*}
\int_{\frac{T}{2}}^TS^k(q_1q_2x, \dfrac{a_1}{q_1}, \dfrac{a_2}{q_2})dx= \frac{(q_1q_2)^k}{2^{\frac{7}{2}k-1}\pi^k}B_k(\frac{a_1}{q_1}, \frac{a_2}{q_2}) \int_{\frac{T}{2}}^Tx^{\frac{k}{4}}dx\\
+O\big((q_1q_2)^{k} T^{1+\frac{k}{4}-\frac{A_0-k}{4(A_0-2)s(K_0)}+\varepsilon}\big).
\end{multline*}

Then we get the evaluation of the integral $\int_1^TS^k(q_1q_2x, \dfrac{a_1}{q_1}, \dfrac{a_2}{q_2})dx$ by summing over integrals over intervals of the form $[2^{-j}T,2^{1-j}T]$, $1\leq j\ll \log T$.


\end{document}